\documentclass[a4paper,12pt]{article}

\usepackage{amssymb, amsthm, amsmath}
\usepackage{fancybox}
\usepackage[dvips]{color}
\usepackage{amsthm}
\usepackage{version}	% required for `\comment' (yatex added)

\newtheorem{definition}{\sc Definition}[section]
\newtheorem{theorem}[definition]{\sc Theorem}
\newtheorem{lemma}[definition]{\sc Lemma}

\newtheorem{proposition}[definition]{\sc Proposition}
\newtheorem{corollary}[definition]{\sc Corollary}
\newtheorem{remark}[definition]{\sc Remark}

\makeatletter
\renewcommand{\theequation}{%
\thesection.\arabic{equation}}
\@addtoreset{equation}{section}

\def\labelenumi{{\rm (\theenumi)}}

\begin{document}

\title{
Global attractors for doubly nonlinear evolution equations with
non-monotone perturbations
}
\author{{\sc Goro Akagi}\footnote{
Department of Machinery and Control Systems, School of Systems
Engineering, Shibaura Institute of Technology, 307 Fukasaku, Minuma-ku,
Saitama-shi, Saitama 337-8570 Japan (e-mail: g-akagi@sic.shibaura-it.ac.jp).
Supported in part by the Shibaura Institute of Technology grant
for Project Research (No.\ 211459 (2006), 211455 (2007)), and the Grant-in-Aid for 
 Young Scientists (B) (No.\ 19740073), Ministry of Education, Culture,
Sports, Science and Technology.
}
}
\date{}
\maketitle

\begin{abstract}
This paper addresses the analysis of dynamical systems
 generated by doubly nonlinear evolution equations governed by
 subdifferential operators with non-monotone perturbations in a
 reflexive Banach space setting. In order to construct global
 attractors, an approach based on the notion of generalized semiflow is
 employed instead of the usual semi-group approach, since solutions of
 the Cauchy problem for the equation might not be unique. Moreover, the
 preceding abstract theory is applied to a generalized Allen-Cahn
 equation as well as a semilinear parabolic equation with a
 nonlinear term involving gradients.\\[5mm]
\end{abstract}

\noindent
%\begin{keyword}
% keywords here, in the form: keyword \sep keyword
%\end{keyword}

\section{Introduction}

Let $V$ and $V^*$ be a real reflexive Banach space and its dual space,
respectively, and let $H$ be a Hilbert space whose dual space $H^*$ is
identified with itself such that
\begin{eqnarray}\label{ps}
 && V \hookrightarrow H \equiv H^* \hookrightarrow V^*
\end{eqnarray}
with continuous and densely defined canonical injections.
Let $\varphi$ and $\psi$ be proper lower semi-continuous convex functions
from $V$ into $(-\infty,\infty]$, and let $\partial_V \varphi, \partial_V
\psi :V \to 2^{V^*}$ be subdifferential operators of $\varphi$ and
$\psi$ respectively. Moreover, let $B$ be a (possibly) non-monotone and
multi-valued operator from $V$ into $V^*$.

We deal with the following Cauchy problem (denoted by (CP)) for a doubly nonlinear
evolution inclusion:
\begin{eqnarray}
&& \partial_V \psi ({u'}(t)) + \partial_V \varphi(u(t)) + \lambda B(u(t)) \ni f \ \mbox{ in }V^*,
\quad 0 < t < \infty,
\label{ee}\\
&& u(0) = u_0,
\label{ic}
\end{eqnarray}
where $f \in V^*$ and $u_0 \in D(\varphi) := \{u \in
V; \varphi(u) < \infty\}$ are given data, with a parameter $\lambda \in
[0, 1]$, which controls the smallness of the perturbation (i.e., no
perturbation if $\lambda = 0$ and no restriction if $\lambda = 1$). 
Equation (\ref{ee}) can be regarded as a perturbation problem for a
monotone system $\partial_V \psi(u'(t)) + \partial_V \varphi(u(t)) \ni
f$ in $V^*$, which has been well studied
(see, e.g.,~\cite{Arai},~\cite{Barbu75},~\cite{Colli},~\cite{CV},~\cite{MT04},~\cite{AFK04,AFK05},~\cite{Roubicek},~\cite{Mielke05},~\cite{RMS08} and the
references given therein), and we emphasize that the operator $B(\cdot)$
is non-monotone and has no gradient structure.
We first review some existence result of global (in time) strong
solutions of (CP) for any $\lambda \in [0, 1]$ (i.e., no restriction on
$\lambda$) under appropriate conditions
such as the coerciveness and the boundedness of $\partial_V \psi$, the
precompactness of sub-level sets of $\varphi$, and the boundedness and
the compactness of $B$ (see Theorem \ref{GE:T:GE}). 

The first purpose of this paper is to study the
long-time behavior of global solutions for (CP), in particular, the
existence of global attractors. There are a number of works for
the study of dynamical systems generated by evolution equations (see,
e.g.,~\cite{Hale},~\cite{LadyAt},~\cite{BaVi},~\cite{Temam}). Most
of them are concerned with single-valued dynamical systems, that is, the
case that the solution is uniquely determined by initial data, and they
are based on a usual semi-group approach.
However, the doubly nonlinear problem (CP) might admit multiple
solutions; indeed, Colli~\cite{Colli} gave a simple example that (CP)
possesses multiple solutions even if $B \equiv 0$ and $\partial_V
\varphi$ is linear. Moreover, the author also exhibited the
non-uniqueness of solutions for some doubly nonlinear problems in~\cite{G09}.
Hence (CP) possibly generates a dynamical system with
a multi-valued evolution operator, which is called {\it multi-valued
dynamical system}. There are several general theories for handling
multi-valued dynamical systems
in~\cite{Sel73},~\cite{ChVi95},~\cite{Ball},~\cite{MeVa98} (see
also~\cite{Sel96},~\cite{MeVa00},~\cite{Ymz02},~\cite{CMR03},~\cite{Ymz04},
~\cite{Ball04},~\cite{CF06},~\cite{Sch07},~\cite{Che09} and references
therein). 
We particularly employ the notion of generalized semiflow proposed by
J.M.~Ball~\cite{Ball} to analyze the long-time behavior of solutions for
(CP).

There are only a few contributions to the long-time behavior of
solutions for doubly nonlinear problems such as (CP). 
In a Hilbert space setting (i.e., $V = V^* = H$), Segatti~\cite{Segatti}
proved the global (in time) solvability of (CP) with $B(u) = -u$ under 
a linear growth condition for $\partial_H \psi$, and moreover, he proved
the existence of global attractors by using the notion of generalized
semiflow and establishing a dissipative estimate for a Lyapunov
functional $J(u) := \varphi(u) - \frac{\lambda}{2} |u|_H^2$. 
Moreover, Schimperna, Segatti and Stefanelli~\cite{SSS} studied the
well-posedness and the long-time behavior of solutions for a concrete
doubly nonlinear gradient system in $L^2(\Omega)$ of a non-convex energy
functional.
 
In our setting, the form of the perturbation has
not been explicitly specified, and moreover, it might have no gradient
structure. Furthermore, we work in a Banach space setting (i.e., $V \neq
V^*$). We use
techniques recently developed by the author in~\cite{G11} to handle
severe nonlinearities arising from the double nonlinearity and the
Banach space setting of (CP), and then, we verify that the set of all
solutions of (CP) forms a generalized semiflow in Ball's sense without
any restrictions on the parameter $\lambda \in [0, 1]$ 
(see Theorem \ref{GF:T:GF}). 
Moreover, we prove the existence of global attractors for the
generalized semiflow in some critical
case by imposing a smallness constraint on $\lambda > 0$ (see Theorem
\ref{At:T:At} and Remark \ref{At:R:crit}). 
We also prove the existence of global
attractors for an arbitrary $\lambda \in [0, 1]$ in subcritical cases (see
Corollary \ref{At:C:At}).

The second purpose of the present paper is to apply the preceding theory
to generalized Allen-Cahn equations proposed by Gurtin~\cite{Gurtin}. We
first investigate the asymptotic behavior of solutions $u =
u(x,t)$ for a quasilinear Allen-Cahn equation of the form
\begin{equation}\label{I:gAC-1}
 \alpha (u_t) - \Delta_m u + \partial_r W(x, u) \ni
  f \quad \mbox{ in } \ \Omega \times (0, \infty),
\end{equation}
where $\Omega$ is a bounded domain of $\mathbb{R}^N$, 
$f: \Omega \to \mathbb{R}$ is given, $\alpha(r) =
|r|^{p-2}r$ with $p \geq 2$ and $\Delta_m$ stands
for the so-called $m$-Laplace operator given by
$$
\Delta_m \phi(x) = \nabla \cdot \left( |\nabla \phi(x)|^{m-2}\nabla \phi(x)
\right),
\quad 1 < m < \infty.
$$ 
Moreover, $W(x,r)$ denotes a double-well potential given by
\begin{equation*}
 W(x,r) = j(r) + \lambda \int^r_0 g(x, \rho) d\rho
  \quad \mbox{ for } \ x \in \Omega, \ r \in \mathbb{R}
\end{equation*}
with a proper lower semicontinuous convex function $j : \mathbb{R} \to
[0, +\infty]$, $\lambda \in [0, 1]$ and a Carath\'eodory function $g : \Omega
\times \mathbb{R} \to \mathbb{R}$ possibly non-monotone in the second
variable. Then $\partial_r W(x,r)$ stands for the derivative in $r$
of the potential $W(x,r)$, more precisely,
$$
\partial_r W(x,r) = \partial j(r) + \lambda g(x,r)
\quad \mbox{ for } \ x \in \Omega \ \mbox{ and } \ r \in \mathbb R,
$$
where $\partial j$ is the subdifferential of $j$.
We reduce the initial-boundary value problem for (\ref{I:gAC-1}) into
(CP) in $V = L^p(\Omega)$ by setting $\partial_V \psi(u) =
\alpha(u(\cdot))$, $\partial_V \varphi(u) = -\Delta_m u +
\partial j(u(\cdot))$ and $B(u) = g(\cdot, u(\cdot))$.
It is crucial in our analysis to set $V = L^p(\Omega)$ so that the
mapping $\partial_V \psi: u \mapsto \alpha(u(\cdot))$ is coercive and
bounded from $V$ into $V^*$.
%because these properties play crucial roles to
%treat (CP).

We next handle the following semilinear problem with a nonlinear term
involving the gradient of $u$:
\begin{equation}\label{I:gAC-2}
\alpha (u_t) - \Delta u + N(x, u, \nabla u) = f
 \quad \mbox{ in } \ \Omega \times (0, \infty),
\end{equation}
where $N(x,r,{\bf p}) = \partial j(r) + \lambda h(x, r, {\bf p})$ with a
Carath\'eodory function $h$ from $\Omega \times \mathbb{R} \times
\mathbb{R}^N$ into $\mathbb{R}$.
This problem could not be reduced into any gradient systems
(cf.~\cite{Segatti},~\cite{SSS},~\cite{RSS}); however, \eqref{I:gAC-2}
can be regarded as a perturbation problem of a gradient system, and it
also falls within our abstract setting.

This paper is composed of six sections. 
In Section \ref{S:P}, we briefly review the theory of generalized
semiflow. In Section \ref{S:GE}, we discuss
the existence of global (in time) solutions for (CP). Then we start our
analysis of the long-time behavior of solutions for (CP). We first
define the set $\mathcal{G}$ of all solutions for (CP) and prove that it
forms a generalized semiflow in a metric space $X$ in Section
\ref{S:GF}. Moreover, we verify the existence of a global attractor for
the generalized semiflow $\mathcal{G}$ in Section \ref{S:At}. 
Applications of the preceding results to
generalized Allen-Cahn equations are given in Section~\ref{S:App}.\\[3mm]
{\bf Notation.} Let $I$ be a section of $\mathbb{R}$ and let $E$ be a
set. We then denote by $AC(I;E)$ (respectively, $AC(I)$) the set of all $E$-valued
(respectively, real-valued) absolutely continuous functions defined on
$I$. Moreover, $C$ stands for a non-negative constant independent of the
elements of the corresponding space and set and may vary from line to
line.

\section{Theory of generalized semiflow}\label{S:P}

The notion of generalized semiflow was first introduced by
J.M.~Ball~\cite{Ball}. He also defined global attractors for generalized
semiflows and provided a criterion of the existence of
global attractors. We first recall the definition of generalized semiflow.
\begin{definition}\label{P:D:GF}
Let $X$ be a metric space with metric $d_X = d_X(\cdot,\cdot)$.
 A family $\mathcal{G}$ of maps $\varphi : [0,\infty) \to X$ is said to be a
 generalized semiflow in $X$, if the following four conditions are all
 satisfied\/{\rm :}
\begin{description}
 \item[(H1)] {\rm (Existence)} for each $x \in X$ there exists $\varphi
	    \in \mathcal{G}$ such that $\varphi(0) = x${\rm ;}
 \item[(H2)] {\rm (Translation invariance)} if $\varphi \in \mathcal{G}$
	    and $\tau \geq 0$, then the map $\varphi^\tau$ also
	    belongs to $\mathcal{G}$,
	    where $\varphi^\tau(t) := \varphi(t + \tau)$ for $t \in [0,
	    \infty)${\rm ;}
 \item[(H3)] {\rm (Concatenation invariance)} if $\varphi_1, \varphi_2 \in
	    \mathcal{G}$ and $\varphi_2(0) = \varphi_1(\tau)$ at some
	    $\tau \geq 0$, then the map $\psi$, 
	    the concatenation of $\varphi_1$ and $\varphi_2$ at $\tau$, 
	    defined by
	    \begin{eqnarray*}
	     \psi(t) := \left\{
			 \begin{array}{ll}
			  \varphi_1(t) \ & \mbox{ if }\ t \in [0,
			   \tau], \\
			  \varphi_2(t - \tau) \ 
			   & \mbox{ if }\ t \in (\tau, \infty)
			 \end{array}
			\right.
	    \end{eqnarray*}
	    also belongs to $\mathcal{G}${\rm ;}
 \item[(H4)] {\rm (Upper semicontinuity)} if $\varphi_n \in
	    \mathcal{G}$, $x \in X$ and $\varphi_n(0) \to x$ in $X$,
	    then there exist a subsequence $(n')$ of $(n)$ 
	    and $\varphi \in \mathcal{G}$ such that $\varphi_{n'} (t)
	    \to \varphi(t)$ for each $t \in [0, \infty)$.
\end{description}
\end{definition}

Let $\mathcal{G}$ be a generalized semiflow in a metric space $X$.
We define a mapping $T(t) : 2^X \to 2^X$ by
\begin{equation}\label{P:T-def}
 T(t) E := \left\{
	    \varphi(t) ;\ \varphi \in \mathcal{G} \ \mbox{ and } \
	    \varphi(0) \in E
	   \right\}
\quad \mbox{ for } \ E \subset X
\end{equation}
for each $t \geq 0$.
Then one can check from (H1)--(H3) that $(T(t))_{t \geq 0}$ satisfies the semi-group
properties, that is, (i) $T(0)$ is the identity mapping in $2^X$; (ii)
$T(t) T(s) = T(t + s)$ for all $t,s \geq 0$. 

Moreover, global attractors for generalized semiflows are defined as
follows.
\begin{definition}
Let $\mathcal{G}$ be a generalized semiflow in a metric space $X$ and
 let $(T(t))_{t \geq 0}$ be the family of mappings defined as in {\rm (\ref{P:T-def})}.
A set $\mathcal{A} \subset X$ is said to be a global attractor for the
 generalized semiflow $\mathcal{G}$ if the following {\rm (i)--(iii)} hold.
\begin{enumerate}
\item[\rm (i)] $\mathcal{A}$ is compact in $X${\rm ;}
\item[\rm (ii)] $\mathcal{A}$ is invariant under $T(t)$, i.e., $T(t)
	     \mathcal{A} = \mathcal{A}$, for all $t \geq 0${\rm ;}
\item[\rm (iii)] $\mathcal{A}$ attracts any bounded subsets $B$ of
	     $X$ by $(T(t))_{t \geq 0}$, i.e., 
	     \begin{equation*}
	      \displaystyle \lim_{t \to
	     \infty}\mathrm{dist}(T(t)B, \mathcal{A}) = 0,
	     \end{equation*}
	     where $\mathrm{dist}(\cdot, \cdot)$ is defined by
	     \begin{equation*}
	      \mathrm{dist}(A,B) := \sup_{a \in A} \inf_{b \in B} d_X(a, b)
	       \quad \mbox{ for } \ A, B \subset X.
	     \end{equation*} 
\end{enumerate}
\end{definition}

As in the standard theory of dynamical systems for semi-group
operators, one can also introduce the notion of $\omega$-limit set.
\begin{definition}
 Let $\mathcal{G}$ be a generalized semiflow in a metric space $X$. For
 $E \subset X$, the $\omega$-limit set of $E$ for $\mathcal{G}$ is given as follows.
\begin{eqnarray*}
 \omega(E) := \Big\{
  x \in X ;\  \mbox{there exist sequences } (\varphi_n)
  \mbox{ in } \mathcal{G} \mbox{ and } (t_n)
  \mbox{ on } [0, \infty) \quad\\ 
  \mbox{such that } \varphi_n(0) \mbox{ is bounded and belongs to } E \mbox{ for } n \in \mathbb{N}, \quad\\
 t_n \to \infty \mbox{ and } \varphi_n(t_n) \to x \mbox{ in } X
  \Big\}.
\end{eqnarray*}
\end{definition}

In order to prove the existence of global attractors for generalized
semiflows, we employ the following theorem due to
J.M.~Ball~\cite{Ball}.
\begin{theorem}[J.M.~Ball~\cite{Ball}]\label{P:T:At}
 A generalized semiflow $\mathcal{G}$ in a metric space $X$ has a global attractor
 $\mathcal{A}$ if and only if the following two conditions are
 satisfied.
\begin{enumerate}
 \item[\rm (i)] $\mathcal{G}$ is point dissipative, that is, 
\begin{equation*}
\exists \mbox{a bounded set } B \subset X,\ 
\forall \varphi \in \mathcal{G},\
\exists \tau = \tau (\varphi) \geq 0,\ 
\forall t \geq \tau,\quad
\varphi(t) \in B;
\end{equation*}
 \item[\rm (ii)] $\mathcal{G}$ is asymptotically compact, that is,
	      for any sequences $(\varphi_n)$ in $\mathcal{G}$ and
	      $(t_n)$ on $[0, \infty)$, if $(\varphi_n(0))$ is
	      bounded in $X$ and $t_n \to \infty$,
	      then $(\varphi_n(t_n))$ is precompact in $X$.
\end{enumerate}
Moreover, $\mathcal{A}$ is the unique global attractor for $\mathcal{G}$ 
and given by
\begin{equation*}
 \mathcal{A} 
  = \bigcup \left\{
	     \omega(B) ; \ B \mbox{ is a bounded subset of } X
	    \right\}
  =: \omega(X).
\end{equation*}
Furthermore, $\mathcal{A}$ is the maximal compact invariant subset of
 $X$ under the family of mappings $(T(t))_{t \geq 0}$.
\end{theorem}

The following proposition gives a sufficient condition for the
asymptotic compactness of generalized semiflows.
\begin{proposition}[J.M.~Ball~\cite{Ball}]\label{P:P:asym-cpt}
 Let $\mathcal{G}$ be a generalized semiflow in a metric space $X$. 
If $\mathcal{G}$ satisfies the following conditions\/{\rm :}
\begin{enumerate}
 \item[\rm (i)] $\mathcal{G}$ is eventually bounded, that is,
 for any bounded set $D \subset X$, there exist $\tau = \tau(D) \geq 0$ and a
  bounded set $B = B(D) \subset X$ such that
\begin{equation*}
 \bigcup_{t \geq \tau} T(t)D \subset B,
\end{equation*}
\item[\rm (ii)] $\mathcal{G}$ is compact, that is, for any sequence
	     $(u_n)$ in $\mathcal{G}$, if $(u_n(0))$ is bounded in
	     $X$, then there exists a subsequence $(n')$ of $(n)$
	     such that $(u_{n'}(t))$ is convergent in $X$ for all $t > 0$,
\end{enumerate}
then $\mathcal{G}$ is asymptotically compact.
\end{proposition}

\section{Global existence}\label{S:GE}

In this section, we discuss the existence of global (in time)
solutions for (CP). 
Let us recall our basic setting again: $V$ and $V^*$ are a real
reflexive Banach space and its dual space, respectively, 
and $H$ is a Hilbert space
satisfying \eqref{ps}. In addition, we always assume that $\varphi$ and
$\psi$ are proper (i.e., $\varphi \not\equiv \infty$, $\psi \not\equiv
\infty$), lower semicontinuous and convex functions from $V$
into $(-\infty, \infty]$ (see also (ii) of Remark
\ref{M:R:g-est-fromEQ}). The subdifferential operator $\partial_V
\varphi : V \to 2^{V^*}$ of $\varphi$ is defined by
$$
\partial_V \varphi(u) := \left\{
\xi \in V^* ; \ \varphi(v) - \varphi(u) \geq
\langle \xi, v - u \rangle \ \mbox{ for all } \ v \in D(\varphi)
\right\},
$$
where $D(\varphi) := \{u \in V ; \ \varphi(u) < \infty\}$ is the
effective domain of $\varphi$, with the domain $D(\partial_V \varphi) :=
\{u \in D(\varphi) ; \ \partial_V \varphi(u) \neq \emptyset\}$. Moreover,
$\partial_V \psi$ is also defined in a similar way. It is well known
that every subdifferential operator is maximal monotone in $V
\times V^*$ (see~\cite{BCP},~\cite{B},~\cite{HB1} for more details).

The existence of solutions for (CP) has already been
studied in more general settings by the author~\cite{G11}. 
Throughout this paper we are concerned with strong solutions for (CP)
given as follows.
\begin{definition}\label{D:sol}
For $T \in (0, \infty)$, a function $u \in AC([0,T];V)$ is said to be a strong
 solution of {\rm (CP)} on $[0,T]$, if the following conditions are satisfied\/{\rm :}
\begin{enumerate}
 \item[\rm (i)] $u(0) = u_0$,
 \item[\rm (ii)] there exists a negligible set $N \subset (0, T)$, i.e.,
	      the Lebesgue measure of $N$ is zero, such that
	      $u(t) \in D(\partial_V \varphi)$ and ${u'}(t) \in  D(\partial_V \psi)$ 
		for all $t \in [0,T] \setminus N$, and moreover, there exist sections
		$\eta(t) \in \partial_V \psi({u'}(t))$, $\xi(t) \in \partial_V \varphi(u(t))$
		and $g(t) \in B(u(t))$ such that 
		\begin{eqnarray}
		 && \eta(t) + \xi(t) + \lambda g(t) = f\ \mbox{ in $V^*$\   for all }t \in [0,T] \setminus N,\label{EQ}
		\end{eqnarray}
 \item[\rm (iii)] $u(t) \in D(\varphi)$ for all $t \in [0, T]$, and 
	      the function $t \mapsto \varphi(u(t))$ is absolutely continuous on $[0,T]$.
\end{enumerate}
Furthermore, for $T \in (0, \infty]$, a function $u \in AC([0,T);V)$ is said to be a strong
 solution of {\rm (CP)} on $[0,T)$, if $u$ is a strong solution of
 {\rm (CP)} on $[0,S]$ for every $S \in (0,T)$.
\end{definition}

In order to discuss the existence of solutions for (CP), 
we give basic assumptions with
parameters $p \in (1,\infty)$, $T > 0$, $\varepsilon > 0$ in the
following. 
Here and thereafter, $p'$ denotes the H\"older conjugate of
$p \in (1, \infty)$, i.e., $p' := p/(p - 1) \in (1, \infty)$, and
moreover, the graph of an operator $A$ is denoted by $A$ again, i.e.,
$[u,\xi] \in A$ means that $u \in D(A)$ and $\xi \in A(u)$.
\begin{description}
 \item[(A1)] There exist constants $C_1 > 0$, $C_2 \geq 0$ such
	    that
	    \begin{eqnarray*}
	     && C_1 |u|_V^p \leq \psi(u) + C_2 \quad
	      \mbox{ for all } u \in D(\psi).
	    \end{eqnarray*}
 \item[(A2)] There exist constants $C_3, C_4 \geq 0$ such that
	    \begin{eqnarray*}
	     && |\eta|_{V^*}^{p'} \leq C_3 \psi(u) + C_4 \quad
	      \mbox{ for all } [u,\eta] \in \partial_V \psi.
	    \end{eqnarray*}
 \item[($\Phi$1)] There exist a reflexive Banach space $X_0$ and a non-decreasing
	    function $\ell_1$ in $\mathbb R$ such that $X_0$ is
	   compactly embedded in $V$ and
	    \begin{eqnarray*}
	     && |u|_{X_0} \leq \ell_1( |u|_H + \varphi(u)) 
	      \quad \mbox{ for all } u \in D(\partial_V \varphi).
	    \end{eqnarray*}
\item[${\bf (B1)}_\varepsilon$] $D(\partial_V \varphi) \subset D(B)$. 
	   There exists a constant
	    $c_\varepsilon \geq 0$ 
	   such that
	    \begin{eqnarray*}
	     && |g|_{V^*}^{p'} \leq \varepsilon |\xi|_{V^*}^\sigma 
	      + c_\varepsilon \left(
			       |\varphi(u)| + |u|_V^p + 1 \right) 
	      \quad \mbox{ with } \sigma := \min\{2,p'\}
	    \end{eqnarray*}
	    for all $u \in D(\partial_V \varphi)$, $g \in B(u)$ and $\xi \in \partial_V \varphi(u)$.
 \item[(B2)] Let $S  \in (0,T]$ and let $(u_n)$ and $(\xi_n)$ be
	    sequences in $C([0,S];V)$ and $L^\sigma(0,S;V^*)$ with
	    $\sigma := \min \{2,p'\}$, respectively, such that $u_n \to u$ strongly in $C([0,S];V)$,
	    $[u_n(t), \xi_n(t)] \in \partial_V \varphi$ for a.e. $t \in (0,S)$, and
	    \begin{eqnarray*}
	     \sup_{t \in [0,S]} |\varphi(u_n(t))| + \int^S_0 |u'_n(t)|_H^p dt
	      + \int^S_0 |\xi_n(t)|_{V^*}^\sigma dt \\
	     \mbox{ is bounded for all } n \in \mathbb{N},
	    \end{eqnarray*}
	    and let $(g_n)$ be a sequence in $L^{p'}(0,S;V^*)$ such
	    that $g_n(t) \in B(u_n(t))$ for a.e. $t \in (0,S)$ and $g_n \to g$ weakly in $L^{p'}(0,S;V^*)$.
	    Then $(g_n)$ is precompact in $L^{p'}(0,S;V^*)$ and $g(t) \in B(u(t))$ for a.e. $t \in (0,S)$.
 \item[(B3)] Let $S \in (0,T]$ and $u \in C([0, S];V) \cap W^{1,p}(0,S;H)$ be such
	    that $\textstyle \sup_{t \in [0,S]} |\varphi(u(t))| < \infty$
	    and suppose that there exists $\xi \in L^{p'}(0,S;V^*)$ such
	    that $\xi(t) \in \partial_V \varphi(u(t))$ for a.e. $t \in (0,S)$. 
	    Then there exists a $V^*$-valued strongly measurable
	    function $g$ such that $g(t) \in B(u(t))$ for a.e. $t \in (0,S)$. 
	    Moreover, the set $B(u)$ is convex for all $u \in D(B)$.
\end{description}

Here let us give a remark on assumptions for the non-monotone
multi-valued operator $B : V \to V^*$.
\begin{remark}
{\rm
\begin{enumerate}
 \item[(i)] Condition ${\rm (B1)}_\varepsilon$ provides some growth
	    condition for $B : V \to V^*$ with a constant 
	    $\varepsilon > 0$. Condition (B2) means that the operator 
	     $\mathcal B : u \mapsto B(u(\cdot))$ is compact and closed 
	    in the sense of
	     multivalued operators. Moreover, we also assume (B3) so
	     that the set $\mathcal B(u)$ lies on a proper
	     Bochner-Lebesgue space for $u$ sufficiently
	     regular in time and it is convex.
 \item[(ii)] In case $B$ is single-valued, one can replace (B2) by
	      simpler conditions. The following condition is
	      much simpler, but somewhat restrictive in view of
	      applications to PDEs.
	      \begin{description}
	       \item[${\rm (B2)}_1$] $B$ is locally Lipschitz continuous
			  from $V$ into $V^*$.
	      \end{description}
	      We also give a milder one.
	      \begin{description}
	       \item[${\rm (B2)}_2$] For each $S \in (0, T)$, the operator 
			  $\mathcal B : u \mapsto B(u(\cdot))$ is
			  continuous and compact from 
			  $L^\infty(0,S;X_0) \cap W^{1,p}(0,S;H)$ into 
			  $L^{p'}(0,S;V^*)$,
	      \end{description}
	      where $X_0$ is the Banach space appeared in ($\Phi$1). 
	      Each of ${\rm (B2)}_1$ and ${\rm (B2)}_2$ implies (B2)
	      under ${\rm (B1)}_\varepsilon$ and ($\Phi$1). Moreover,
	      ${\rm (B2)}_2$ can be also derived from a pointwise condition:
	      \begin{description}
	       \item[${\rm (B2)}_3$] There exists a Banach space $W$
			  such that $X_0$ is compactly embedded in $W$ and
			  $W$ is continuously embedded in $H$. Moreover,
			  it holds that
			  $$
			  \left|B(u) - B(v)\right|_{V^*}
			  \leq \ell_2 ( |u|_{X_0} )
			  \ell_3 (|u - v|_W)
			  \quad \mbox{ for all } \ u, v \in X_0
			  $$
			  with non-decreasing functions $\ell_2,\ell_3$
			  in $\mathbb R$ satisfying 
			  $\ell_3(s) \to 0$ as $s \to +0$.
	      \end{description}
	      Indeed, if $(u_n)$ is bounded in 
	      $L^\infty(0,S;X_0) \cap W^{1,p}(0,S;H)$, then by Aubin-Lions's
	      compactness theorem, up to a subsequence, one can ensure
	      that $u_n \to u$ strongly in $C([0,S];W)$. Thus
	      $B(u_n(\cdot)) \to B(u(\cdot))$ strongly in
	      $L^\infty(0,S;V^*)$.
 \item[(iii)] If $B$ is single-valued and either ${\rm (B2)}_1$ or
	      ${\rm (B2)}_2$ holds, then (B3) is automatically verified.
\end{enumerate}
}
\end{remark}

\begin{remark}
{\rm
Assumption (A2) implies that
\begin{description}
 \item[\bf $\mbox{(A2)}'$] There exists a constant
	    $C_5 \geq 0$ such that
	    \begin{equation*}
	     \frac 1 2 \psi(u) \leq \langle \eta, u \rangle + C_5
	      \quad \mbox{ for all } [u, \eta]  \in \partial_V \psi.
	    \end{equation*} 
\end{description}
For its proof, we refer the reader to Proposition 3.2 of~\cite{G11}.
}
\end{remark}

Now, applying the abstract theory developed in~\cite{G11}, we can assure
that
\begin{theorem}[Global existence]\label{GE:T:GE}
 Let $p \in (1,\infty)$, $T > 0$ and $\lambda \in [0,1]$ be fixed.
 Suppose that {\rm (A1), (A2), ${\rm (\Phi 1)}$, ${\rm (B1)}_\varepsilon$--(B3)} 
are all satisfied with a sufficiently small $\varepsilon > 0$
{\rm (}the smallness of $\varepsilon$ is determined only by $p$, $C_1$,
 $C_3$ and $C_H := \sup \{|u|_H/|u|_V ;~u \in V \setminus \{0\}\} > 0${\rm )}. 
Then, for all $f \in V^*$ and $u_0 \in D(\varphi)$, there exists at least one strong
 solution $u \in W^{1,p}(0,T;V)$ on $[0,T]$ such that
\begin{equation*}
 \eta, \xi, g \in L^{p'}(0,T; V^*),
\end{equation*}
where $\eta(\cdot)$, $\xi(\cdot)$ and $g(\cdot)$ are the sections of 
$\partial_V \psi(u'(\cdot))$, $\partial_V \varphi(u(\cdot))$
and $B(u(\cdot))$, respectively, as in {\rm (\ref{EQ})}.
\end{theorem}

We close this section with the following remark.
\begin{remark}\label{M:R:g-est-fromEQ}
 {\rm
\begin{enumerate}
 \item[(i)] 
As in \S 6.1 of~\cite{G11}, we have:
\begin{proposition}\label{GE:P:B-bdd}
Let $p \in (1, \infty)$ and $T > 0$  be given. 
Assume that {\rm (A2)} and ${\rm (B1)}_\varepsilon$ hold with an
 arbitrary constant $\varepsilon \in (0, 4^{1-p'})$.
Let $\lambda \in [0, 1]$ and let $u$ be a strong solution of {\rm (CP)}
 on $[0, T]$ and let $g(\cdot)$ denote the section of $B(u(\cdot))$ as in {\rm
 (\ref{EQ})}. Then there exist constants 
$\gamma_\varepsilon = \gamma_\varepsilon(p, C_3) \geq 0$ and
 $M_\varepsilon = M_\varepsilon(p, c_\varepsilon) \geq 0$ such that 
 \begin{eqnarray*}
  && |g(t)|_{V^*}^{p'}   
   \leq M_\varepsilon \left(
	   |f|_{V^*}^{p'} + C_4 
	   + |\varphi(u(t))| + |u(t)|_V^p
	   + 1	      \right)
   + \gamma_\varepsilon \psi \left(
		  u'(t)
		 \right) 
 \end{eqnarray*}
for a.e. $t \in (0, T)$. Here $\gamma_\varepsilon$ depends on
 $\varepsilon$ as follows\/{\rm :}
\begin{equation*}
 \gamma_\varepsilon \to 0 \quad \mbox{ as } \ \varepsilon \to +0.
\end{equation*}
\end{proposition}
\begin{proof}
By ${\rm (B1)}_\varepsilon$, we have
\begin{eqnarray*}
 |g(t)|_{V^*}^{p'} 
  &\leq& \varepsilon |\xi(t)|_{V^*}^\sigma 
  + c_\varepsilon \left(
		   |\varphi(u(t))| + |u(t)|_V^p + 1
		  \right)
\\
 &\leq&  4^{p'- 1} \varepsilon \left(
				|f|_{V^*}^{p'} + |\eta(t)|_{V^*}^{p'} + |g(t)|_{V^*}^{p'} + 1
			       \right)
 \\
 && + c_\varepsilon \left(
		   |\varphi(u(t))| + |u(t)|_V^p + 1
		  \right),
\end{eqnarray*}
which together with (A2) implies
\begin{eqnarray*}
\left(
1 - 4^{p' - 1} \varepsilon 
\right) |g(t)|_{V^*}^{p'} 
 &\leq& 4^{p'- 1} \varepsilon \left(
			       |f|_{V^*}^{p'} + C_4 + 1
			       \right)
 + 4^{p'- 1} \varepsilon C_3 \psi(u'(t))
\\
&& \quad + c_\varepsilon \left(
		   |\varphi(u(t))| + |u(t)|_V^p + 1
		  \right).
\end{eqnarray*} 
By setting $\gamma_\varepsilon := 4^{p' - 1} \varepsilon C_3/ (1 - 4^{p'-1} \varepsilon) > 0$, 
we can derive our desired result.
\end{proof}
\item[(ii)] We can assume $\varphi \geq 0$ and $\psi \geq 0$ without
	     any loss of generality (see Remark 3.8
	     of~\cite{G11}). In the rest of this paper, we always assume so.
\end{enumerate}
}
\end{remark}

\section{Formation of a generalized semiflow}\label{S:GF}

Our analysis of the large-time behavior of solutions for (CP) is
based on the theory of generalized semiflow briefly reviewed in \S
\ref{S:P}. Let $X := D(\varphi)$ be a metric space equipped with the distance
$d_X(\cdot , \cdot)$ defined by
\begin{equation*}
d_X(u, v) := |u - v|_V + |\varphi(u) - \varphi(v)|
  \quad \mbox{ for }\ u, v \in X
\end{equation*}
and set
\begin{eqnarray*}
\mathcal{G} := \{u \in AC([0, \infty); X) ;\
  u \mbox{ is a strong solution of (CP) on } [0, \infty) \quad \\
  \mbox{with some } u_0 \in X\}.
\end{eqnarray*}
Then our result of this section reads,
\begin{theorem}\label{GF:T:GF}
Let $p \in (1, \infty)$ and $\lambda \in [0,1]$ be fixed. 
Assume that {\rm (A1), (A2), ($\Phi$1), ${\rm (B1)}_\varepsilon$--(B3)}
 are satisfied with a sufficiently small $\varepsilon > 0$ for any $T > 0$ 
{\rm (}the smallness of $\varepsilon$ is determined only by $p$, $C_1$,
 $C_3$ and $C_H${\rm )}.
 Then $\mathcal{G}$ is a generalized semiflow in $X$.
\end{theorem}

\begin{proof}
It suffices to check the four conditions (H1)--(H4) (see Definition
 \ref{P:D:GF}). Since (H2) and (H3) follow immediately, we give proofs
 only for (H1), (H4).\\[3mm]
\noindent
{\it Proof of {\rm (H1)}.} 
Let $u_0 \in X$ be given. Then due to Theorem \ref{GE:T:GE}, for any $T
 > 0$, the Cauchy problem (CP) admits at least one strong solution $u_1
 \in W^{1,p}(0, T; V)$. Moreover, since $u_1(T)$ also belongs to
 $D(\varphi)$, we can also obtain a strong solution $u_2 \in W^{1,p}(0, T;
 V)$ of (CP) with $u_0$ replaced by $u_1(T)$. Iterating the argument
 above, we can construct a sequence $(u_k)_{k \in \mathbb{N}}$ in $W^{1,p}(0,T;V)$
 such that $u_{k+1}$ becomes a strong solution of (CP) with the initial
 condition $u_{k+1}(0) = u_k(T)$ for each $k \in \mathbb{N}$.

Furthermore, define $u : [0,\infty) \to V$ by
 \begin{eqnarray*}
  && u(t) := u_k(t - (k - 1)T ) \quad \mbox{ if }\ t \in [(k - 1)T, kT].
 \end{eqnarray*}
Then $u \in W^{1,p}_{loc}([0, \infty) ; V)$ and $\varphi(u(\cdot)) \in AC([0, \infty))$,
 and moreover, the restriction of $u$ on $[0, S]$
 becomes a strong solution of (CP) for any $S > 0$. 
Therefore $u$ is a strong solution of (CP) on $[0, \infty)$, which
 implies that $u \in \mathcal{G}$ and $u(0) = u_0$. Thus (H1) follows.\\[3mm]
\noindent
{\it Proof of {\rm (H4)}.}
Let $u_n \in \mathcal{G}$ and $v \in X$ be such that $u_n(0) \to v$ in
 $X$. Then multiplying (\ref{ee}) with $u = u_n$ by $u_n'(t)$, we have
\begin{eqnarray*}
 && \left\langle \eta_n(t), u_n'(t) \right\rangle
  + \left\langle \xi_n(t),  u_n'(t) \right\rangle
  + \lambda \left\langle g_n(t), u_n'(t) \right\rangle
  = \left\langle f, u_n'(t) \right\rangle
\end{eqnarray*}
with sections $\eta_n(t) \in \partial_V \psi (u_n'(t))$, 
$\xi_n(t) \in \partial_V \varphi(u_n(t))$, $g_n(t) \in B(u_n(t))$
for a.e. $t \in (0, \infty)$. 
We then derive from ${\rm (A2)}'$ and the
 chain rule for subdifferentials that
\begin{eqnarray*}
\frac 1 2 \psi \left(
		      u_n'(t)
		 \right) 
 - C_5
 + \frac{d}{dt} \varphi(u_n(t))
&=& \left\langle f - \lambda g_n(t),  u_n'(t) \right\rangle
\\
&\leq& \left(
      |f|_{V^*} + \lambda |g_n(t)|_{V^*}
     \right) 
|u_n'(t)|_{V}.
\end{eqnarray*}
Furthermore, by Young's inequality and (A1),
\begin{eqnarray}\label{GF:xdu-p1}
\lefteqn{
\frac 1 4 \psi \left(
			      u_n'(t)
			     \right)
 + \frac{d}{dt} \varphi(u_n(t))
}\nonumber \\
&\leq& M_1 \left(
	|f|_{V^*}^{p'} + \lambda^{p'} |g_n(t)|_{V^*}^{p'} + C_2 + C_5 + 1
       \right)
\end{eqnarray}
with a constant $M_1$ depending only on $p$ and $C_1$.
By Proposition \ref{GE:P:B-bdd}, we can choose $\varepsilon > 0$ so
 small that $\gamma_\varepsilon = \gamma_\varepsilon (p, C_3) \leq 
(8 M_1)^{-1}$
(hence, the smallness of $\varepsilon$ is determined by $p$, $C_1$
 and $C_3$), and then,
\begin{eqnarray}
|g_n(t)|_{V^*}^{p'}
 &\leq& M_\varepsilon \left(
		      |f|_{V^*}^{p'} + C_4 + \varphi(u_n(t)) + |u_n(t)|_V^p + 1
		     \right)
 \nonumber \\
 && + \frac{1}{8 M_1} \psi \left(
		 u_n'(t)
		\right) 
\label{GF:B1}
\end{eqnarray}
for a.e. $t \in (0, \infty)$. Then it follows from (\ref{GF:xdu-p1}) and (\ref{GF:B1}) that
\begin{eqnarray}
 \frac 1 {16} \psi(u_n'(t))
  + \lefteqn{
  \frac{C_1}{16} |u_n'(t)|_V^p
  + \frac{d}{dt} \varphi(u_n(t))
  } \nonumber\\
 &\leq& 
  C \left(
     |f|_{V^*}^{p'} + 1
    \right)
  + \lambda^{p'} M_2 \Big(
	       \varphi(u_n(t)) + |u_n(t)|_V^p
	      \Big),
\label{GF:xdu}
\end{eqnarray}
where $M_2 := M_1 M_\varepsilon$ is independent of $\lambda$ (this fact will be used
 in \S \ref{S:At}). 

Here, for each $\delta > 0$,
 we can take a constant $C_{\delta,p}$ depending only on $\delta$ and
 $p$ such that
\begin{eqnarray*}
 &&
  \frac{d}{dt} |u_n(t)|_V^p 
  \leq \delta \left|
	       u_n'(t)
	      \right|_V^p
  + C_{\delta,p} |u_n(t)|_V^p.
\end{eqnarray*}
Put $\delta = C_1 / 16$. Therefore
\begin{eqnarray}
\lefteqn{
  \frac 1{16} \psi(u_n'(t))
  + \frac{d}{dt} 
  \Big(
   \varphi(u_n(t)) + |u_n(t)|_V^p
  \Big)
} \nonumber \\
 &\leq& C \left(
	|f|_{V^*}^{p'} + 1
       \right)
  + C \Big(
       \varphi(u_n(t)) + |u_n(t)|_V^p
      \Big)
\label{GF:xdu2}
\end{eqnarray}
for a.e. $t \in (0, \infty)$. 
Hence integrate both sides over $(0, t)$ to get
\begin{eqnarray*}
\lefteqn{
\frac 1{16} \int^t_0 \psi(u_n'(\tau)) d\tau
+
 \varphi(u_n(t)) + |u_n(t)|_V^p
 }\\
 &\leq& 
 \varphi(u_0) + |u_0|_V^p 
 + Ct \left(
       |f|_{V^*}^{p'} + 1
      \right)
\\
&& + C  \int^t_0 \Big(
		  \varphi(u_n(\tau)) + |u_n(\tau)|_V^p
		 \Big) d\tau
\end{eqnarray*}
for all $t \geq 0$. Thus Gronwall's inequality yields
\begin{eqnarray}
\lefteqn{
 \sup_{t \in [0, T]} \Big(
		      \varphi(u_n(t)) + |u_n(t)|_V^p
		     \Big)
}\nonumber \\
 &\leq& \left\{
       \varphi(u_0) + |u_0|_V^p + CT \left(
				       |f|_{V^*}^{p'} + 1      
				     \right)
      \right\} e^{CT}
\label{GF:e:u+phi}
\end{eqnarray}
with an arbitrary positive number $T > 0$.
Hence ($\Phi$1) implies that $(u_n(t))_{n \in \mathbb{N}}$
 is precompact in $V$ for each $t \in [0, T]$. Furthermore, we can also deduce that
\begin{eqnarray}
 && \int^T_0 \psi(u_n'(t)) dt \leq C_T, \label{GF:e:psi}\\
 && \int^T_0 |u_n'(t)|_V^p dt \leq C_T \label{GF:e:du}
\end{eqnarray}
with a constant $C_T \geq 0$. Here and henceforth, $C_T$ denotes a
 non-negative constant
 independent of $n$ and $t$ but possibly depending on $T$ and may vary
 from line to line. 

We particularly put $T = 1$ in (\ref{GF:e:u+phi}) and
 (\ref{GF:e:du}). Then by virtue of Ascoli's compactness theorem, $(u_n)$
 becomes precompact in $C([0, 1]; V)$. Hence there exist a subsequence 
$(n^{(1)}_k)$ of $(n)$ and
 a function $u^{(1)} \in C([0, 1];V)$ such that
\begin{eqnarray*}
 && u_{n^{(1)}_k} \to u^{(1)} \quad \mbox{ strongly in } C([0, 1];V).
\end{eqnarray*}
Moreover, from the arbitrariness of $T$ in (\ref{GF:e:u+phi}) and
 (\ref{GF:e:du}), we can iteratively take a function $u^{(m+1)} \in C([0,
 m+1] ; V)$ and a subsequence $(n^{(m+1)}_k)$ of $(n^{(m)}_k)$ such that
\begin{eqnarray*}
 && u_{n^{(m+1)}_k} \to u^{(m+1)} \quad \mbox{ strongly in } C([0, m+1];V)
\end{eqnarray*}
for each $m \in \mathbb{N}$.
Here we remark that the restriction of $u^{(m+1)}$ on $[0, m]$ coincides
 with $u^{(m)}$ for every $m \in \mathbb{N}$. Furthermore, let $(n_k)$
 be the diagonal part of the double sequence $(n^{(m)}_k)_{k,m \in \mathbb{N}}$,
that is, $n_k := n^{(k)}_k$. Then we obtain
\begin{eqnarray}\label{GF:c:ut}
 && u_{n_k}(t) \to u(t) \quad \mbox{ strongly in } V \ \mbox{ for each } \ t \in [0, \infty), 
\end{eqnarray}
where $u \in C([0, \infty);V)$ is given by 
\begin{eqnarray*}
&& u(t) = u^{(m)}(t) \quad \mbox{ if } \
t \in [0, m] \ \mbox { for some } \ m \in \mathbb{N}. 
\end{eqnarray*}

We next prove that $u \in \mathcal{G}$, that is, $u$ is a strong
 solution of (CP) on $[0, \infty)$. Recall
 (\ref{GF:e:u+phi})--(\ref{GF:e:du}) with an arbitrary $T > 0$ and use
 (A2), (\ref{GF:B1}) and (\ref{EQ}) to get
\begin{eqnarray}
 && \int^T_0 |\eta_n(t)|_{V^*}^{p'} dt \leq C_T, \label{GF:e:dpsi}\\
 && \int^T_0 |g_n(t)|_{V^*}^{p'} dt \leq C_T, \label{GF:e:B}\\
 && \int^T_0 |\xi_n(t)|_{V^*}^{p'} dt \leq C_T \label{GF:e:dphi} 
\end{eqnarray}
for all $n \in \mathbb{N}$. Hence we can take a subsequence of
 $(n_k)$, which will be denoted by the same letter, such that
\begin{eqnarray}
 u_{n_k} \to u && \mbox{ weakly in } W^{1,p}(0,T;V),
  \label{GF:c:du}\\
 && \mbox{ strongly in } C([0,T];V),
  \label{GF:c:u}\\
 \eta_{n_k} \to \eta && \mbox{ weakly in } L^{p'}(0,T;V^*),
  \label{GF:c:dpsi}\\
 g_{n_k} \to g && \mbox{ weakly in } L^{p'}(0,T;V^*),
  \label{GF:c:B}\\
 \xi_{n_k} \to \xi && \mbox{ weakly in } L^{p'}(0,T;V^*)
  \label{GF:c:dphi}
\end{eqnarray}
with some $\eta, g, \xi \in L^{p'}(0,T;V^*)$. 
Therefore from the
 demiclosedness of subdifferential operators we can derive that 
$\xi(t) \in \partial_V \varphi(u(t))$ for a.e. $t \in (0, T)$, and
 moreover, (B2) implies that $g(t) \in B(u(t))$ for a.e. $t \in (0, T)$
 and
\begin{eqnarray}
 g_{n_k} \to g && \mbox{ strongly in } L^{p'}(0,T;V^*).
\label{GF:c:B-s}
\end{eqnarray}
Furthermore, multiply $\eta_{n_k}(t)$ by $u_{n_k}'(t)$ and integrate
 this over $(0, T)$. It then follows that
\begin{eqnarray*}
\lefteqn{
\int^T_0 \left\langle
	      \eta_{n_k}(t), u_{n_k}'(t)
	      \right\rangle dt
}\\
&& = - \varphi(u_{n_k}(T)) + \varphi(u_{n_k}(0))
 + \int^T_0 \left\langle f - \lambda g_{n_k}(t), u_{n_k}'(t) 
	    \right\rangle dt.
\end{eqnarray*}
Therefore we get, by (\ref{GF:c:du}), (\ref{GF:c:u}) and (\ref{GF:c:B-s}),
\begin{eqnarray}
\lefteqn{
\limsup_{n_k \to \infty}
\int^T_0 \left\langle
	      \eta_{n_k}(t), u_{n_k}'(t)
	      \right\rangle dt
}\nonumber \\
&\leq& - \varphi(u(T)) + \varphi(v)
 + \int^T_0 \left\langle f - \lambda g(t), u'(t) \right\rangle dt.
\label{GF:lims_dphidu}
\end{eqnarray}
Here we also used the fact that $u_n (0) \to v$ in $X$, in particular,
 $\varphi(u_n(0)) \to \varphi(v)$. Thus by Proposition 1.1
 of~\cite[Chap.~II]{B} (see also~\cite{BCP}), 
we obtain $\eta(t) \in \partial_V \psi(u'(t))$ for
 a.e. $t \in (0, T)$, which implies that $u$ is a strong solution of
 (CP) with $u_0 = v$ on $[0, T]$. From the arbitrariness of $T$, we can
 also verify that $u$ is a strong solution of (CP) on $[0, \infty)$. 
Hence $u$ belongs to $\mathcal{G}$.

Finally, we show that $\varphi(u_{n_k}(t)) \to \varphi(u(t))$ for every
 $t \in [0, \infty)$, by taking a subsequence of $(n_k)$ independent of
 $t$.
To do so, we first derive the convergence for a.e. $t \in (0, \infty)$.
\begin{lemma}\label{GF:L:phi-conv}
It follows that
\begin{eqnarray*}
 && \liminf_{n_k \to \infty} \varphi(u_{n_k}(t)) = \varphi(u(t))
  \quad \mbox{ for a.e. } \ t \in (0, \infty).
\end{eqnarray*}
\end{lemma}

\begin{proof}
Let $T > 0$ be arbitrarily given.
By the definition of subdifferentials, we see
\begin{equation*}
\varphi(u_n(t)) 
 \leq
  \varphi(u(t)) + |\xi_n(t)|_{V^*} |u_n(t) - u(t)|_V.
\end{equation*}
Now, by (\ref{GF:e:dphi}), Fatou's lemma ensures that
\begin{eqnarray*}
 && p(\cdot) := \liminf_{n \to \infty} |\xi_n(\cdot)|_{V^*}^{p'} \in L^1(0, T).
\end{eqnarray*}
Hence $p(t) < \infty$ for a.e. $t \in (0, T)$. Thus we get, by (\ref{GF:c:ut}),
\begin{eqnarray*}
\liminf_{n_k \to \infty} \varphi(u_{n_k}(t)) 
 &\leq& 
 \varphi(u(t)) 
 + \liminf_{n_k \to \infty}
 \left(
  |\xi_{n_k}(t)|_{V^*} |u_{n_k}(t) - u(t)|_V
 \right)
 \\
 &\leq&
 \varphi(u(t)) 
 + p(t) \left(
	 \lim_{n_k \to \infty} |u_{n_k}(t) - u(t)|_V
	\right)\\
 &=&
  \varphi(u(t))
\quad \mbox{ for a.e. } t \in (0, T).
\end{eqnarray*}
Combining this fact with the lower semicontinuity of $\varphi$, we can prove this lemma.
\end{proof}

\noindent
{\it Continuation of proof of {\rm (H4)}}.
We next exhibit the convergence of $\varphi(u_{n_k}(t))$ at every $t \in [0, \infty)$.
Recalling (\ref{GF:xdu}), we find
\begin{eqnarray*}
 && \frac{d\zeta_n}{dt} (t) \leq 0 \quad \mbox{ for a.e. } t \in (0, \infty),
\end{eqnarray*}
where $\zeta_n$ is an absolutely continuous function from $[0, \infty)$ into $\mathbb{R}$ given by
\begin{eqnarray*}
 && \zeta_n(t) := \varphi(u_n(t)) - C t \left(
					 |f|_{V^*}^{p'} + 1
					\right)
 - \lambda^{p'} M_2 \int^t_0 \Big(
	       \varphi(u_n(\tau)) + |u_n(\tau)|_V^p
	      \Big) d\tau
\end{eqnarray*}
for $t \in [0, \infty)$.
Hence $\zeta_n$ is non-increasing on $[0, \infty)$.
We further exploit the following lemma (see Lemma 3.3.3
 of~\cite{GradFlow} for more details).

\begin{lemma}[Helly]\label{GF:L:z-conv}
If $f_n$ is a non-increasing function on $[0, \infty)$, then 
there exist a non-increasing function $g : [0, \infty) \to [-\infty, \infty]$ 
and a subsequence of $(n')$ of $(n)$ independent of $t$ such that
\begin{eqnarray*}
 && \lim_{n' \to \infty} f_{n'}(t) = g (t)
  \quad \mbox{ for all } \ t \in [0, \infty).
\end{eqnarray*}
\end{lemma}

\noindent
{\it Continuation of proof of {\rm (H4)}.}
We apply Lemma \ref{GF:L:z-conv} to our situation with $f_n = \zeta_n$
 and write $\phi$ for the non-increasing function $g$. Then
\begin{eqnarray*}
 && 
  \lim_{n_k \to \infty} \zeta_{n_k}(t) 
  = \phi(t)
  \quad \mbox{ for all } \ t \in [0, \infty).
\end{eqnarray*}
It remains only to reveal the representation of $\phi$.

\begin{lemma}\label{GF:L:phi}
For each $t \in [0, \infty)$, it follows that
 \begin{eqnarray*}
  && \phi(t) = \varphi(u(t)) - C t \left(
				    |f|_{V^*}^{p'} + 1
				   \right)
  - \lambda^{p'} M_2 \int^t_0 \Big(
		\varphi(u(\tau)) + |u(\tau)|_V^p
	       \Big) d\tau.
 \end{eqnarray*}
\end{lemma}

\begin{proof}
Let $T > 0$ be fixed. We can then derive that
\begin{eqnarray*}
 \int^t_0 \varphi(u_{n_k}(\tau)) d\tau \to \int^t_0 \varphi(u(\tau)) d\tau
\quad \mbox{ for all }\ t \in [0, T]
\end{eqnarray*}
from the definition of subdifferentials together with (\ref{GF:c:u}) and
 (\ref{GF:c:dphi}). Hence the definition of $\zeta_n$ and Lemma
 \ref{GF:L:phi-conv} yield
\begin{eqnarray*}
\liminf_{n_k \to \infty} \zeta_{n_k}(t)
= \zeta(t) \quad \mbox{ for a.e. }\ t \in (0, T),
\end{eqnarray*}
where $\zeta \in AC([0, \infty))$ is given by
\begin{equation*}
\zeta(t)
:= \varphi(u(t)) - C t \left(
				    |f|_{V^*}^{p'} + 1
				   \right)
  - \lambda^{p'} M_2 \int^t_0 \Big(
		\varphi(u(\tau)) + |u(\tau)|_V^p
	       \Big) d\tau.
\end{equation*}
Therefore we can obtain $\phi(t) = \zeta(t)$ for a.e. $t \in (0,
 \infty)$ from the arbitrariness of $T > 0$. 

Furthermore, since $\zeta$ is continuous and $\phi$ is non-increasing on
 $[0,\infty)$, we claim that
\begin{eqnarray*}
 && \phi(t) = \zeta(t) \quad \mbox{ for all }\ t \in [0, \infty).
\end{eqnarray*}
Indeed, in case $t \in (0, \infty)$, since $\zeta = \phi$ a.e. in
 $(0, \infty)$, for each $n \in \mathbb{N}$ there exist $t_n \in (t -
 1/(2n), t) \cap [0, \infty)$ and $s_n \in (t, t + 1/(2n))$ such that $\phi = \zeta$ at
 $t_n$ and $s_n$. Since $\phi$ is non-increasing, we find
 that $\phi(t_n) \geq \phi(t) \geq \phi(s_n)$. Hence letting
 $n \to \infty$, we deduce that $\phi(t) = \zeta(t)$, by the continuity
 of $\zeta$. In case $t = 0$, it follows that $\phi(0) = \lim_{n_k \to \infty} \zeta_{n_k}(0) 
= \lim_{n_k \to \infty} \varphi (u_{n_k}(0)) = \varphi(v) = \zeta(0)$
 from the fact that $\varphi(u_n(0)) \to \varphi(v)$.
\end{proof}

\noindent
{\it Continuation of proof of {\rm (H4)}.}
From Lemma \ref{GF:L:phi}, we conclude that $u_{n_k}(t) \to u(t)$ in $X$
 for each $t \in [0, \infty)$. Thus (H4) follows.
\end{proof}

\section{Existence of global attractors}\label{S:At}

In this section, we prove the existence of global attractors for the
generalized semiflow $\mathcal{G}$. To this end, we introduce the
following structure condition, which plays an essential role in our
proof (see also Remark \ref{At:R:crit}).
\begin{description}
 \item[(S1)] For all $\lambda \in [0, \bar \lambda]$ with some 
	    $\bar \lambda \in (0, 1]$, there exist constants $\alpha > 0$
	    independent of $\lambda$ and  $C_6 = C_6(\lambda) \geq 0$ (possibly depending on
	    $\lambda$) such that
	    \begin{eqnarray*}
	     && \alpha \Big(
			\varphi(u) + |u|_V^p
			\Big)
	     \leq \langle \xi + \lambda g, u \rangle + C_6
	    \end{eqnarray*}
	    for all $u \in D(\partial_V \varphi) \cap D(B)$, 
	    $\xi \in \partial_V \varphi(u)$ and $g \in B(u)$.
\end{description}

Our result reads,

\begin{theorem}\label{At:T:At}
Let $p \in (1, \infty)$ be fixed and assume that {\rm (A1), (A2),
 ($\Phi$1)} and {\rm ${\rm (B1)}_\varepsilon$--(B3)} are satisfied with
 a sufficiently small $\varepsilon > 0$ {\rm (}the smallness of $\varepsilon$
 is determined only by $p$, $C_1$, $C_3$ and $C_H${\rm )} for any $T > 0$
 such that every constant and function appeared in the assumptions are
 independent of $\lambda$. Moreover, suppose that {\rm (S1)} is satisfied. 
Then there exists $\lambda_0 \in (0, \bar \lambda]$ such that if 
 $\lambda \in [0,\lambda_0]$ then
the generalized semiflow $\mathcal{G}$ has a unique global attractor
 $\mathcal{A}$, which is given by
\begin{eqnarray*}
 && \mathcal{A} := \bigcup \left\{
			    \omega (B) ; B \mbox{ is a bounded subset of } X
			   \right\}
 = \omega(X).
\end{eqnarray*}
Furthermore, $\mathcal{A}$ is the maximal compact invariant subset of $X$.
\end{theorem}

In order to prove this theorem, thanks to Theorem \ref{P:T:At} and
Proposition \ref{P:P:asym-cpt}, it suffices to prove that the
generalized semiflow $\mathcal{G}$ is point dissipative, eventually
bounded and compact. To do so, we first establish a dissipative estimate
for $u \in \mathcal{G}$ in $X$.

\begin{lemma}\label{At:L:phi-bdd}
Under the same assumptions as in Theorem {\rm \ref{At:T:At}}, 
there exist constants $\lambda_0 \in (0, \bar \lambda]$, $R \geq 0$ and an
 increasing function $T_0(\cdot)$ on $[0, \infty)$ such that 
if $\lambda \in [0, \lambda_0]$, then
\begin{eqnarray}\label{At:phi-bdd}
\varphi(u(t)) + |u(t)|_V^p \leq R
  && \mbox{ for all } \ u_0 \in X, \ 
  u \in \mathcal{G} \mbox{ satisfying } u(0) = u_0
  \nonumber \\
 && \quad \mbox{and } \ t \geq T_0(\varphi(u_0) + |u_0|_V^p).
\end{eqnarray}
\end{lemma}

\begin{proof}
 Let $u_0 \in X$ and let $u \in \mathcal{G}$ be such that $u(0) = u_0$.
By (S1) and (\ref{EQ}), we find that
\begin{eqnarray*}
\alpha \Big(
	    \varphi(u(t)) + |u(t)|_V^p
	    \Big)
 &\leq& 
 \left\langle \xi(t) + \lambda g(t), u(t) \right\rangle + C_6
 \\
 &=&
  \left\langle f - \eta(t), u(t) \right\rangle + C_6
  \\
 &\leq&
  C \left(
   |f|_{V^*}^{p'} + |\eta(t)|_{V^*}^{p'} + 1
  \right)
  + \frac{\alpha}{2} |u(t)|_V^p
\end{eqnarray*}
with some sections $\xi(t) \in \partial_V \varphi(u(t))$, $g(t) \in B(u(t))$,
 $\eta(t) \in \partial_V \psi(u'(t))$
for a.e. $t \in (0, \infty)$. Thus by (A2) it follows that
\begin{eqnarray}
 \frac{\alpha}{2} 
  \Big(
   \varphi(u(t)) + |u(t)|_V^p
  \Big)
 &\leq& 
  C \left(
   |f|_{V^*}^{p'} + 1
  \right)
  + C_7 \psi(u'(t))
\label{At:xu}
\end{eqnarray}
for a.e. $t \in (0, \infty)$ with some constant 
$C_7 = C_7(p,\alpha,C_3) \geq 0$ independent of $\lambda$. 

Recall (\ref{GF:xdu}) with $u_n$ replaced by $u$ and note that
\begin{equation*}
 \frac{C_1}{16} |u'(t)|_V^p 
  \geq \sigma |u'(t)|_V^p
  \geq \sigma \frac{d}{dt} |u(t)|_V^p - \frac{p}{p'} \sigma |u(t)|_V^p,
\end{equation*}
where
\begin{equation*}
 \sigma := \min \left\{
		 \frac{C_1}{16}, 
		 \frac{\alpha p'}{64 p C_7}, 1
		\right\} \in (0, 1].
\end{equation*}
Then it follows that
\begin{eqnarray}
\lefteqn{
 \frac 1{16} \psi(u'(t)) 
+ \frac{d}{dt} \Big(
		\varphi(u(t)) + \sigma |u(t)|_V^p
	       \Big)
}\nonumber \\
 &\leq&
  C \left(
     |f|_{V^*}^{p'} + 1
    \right)
  + \left(
     \lambda^{p'} M_2 + \frac{p}{p'} \sigma
    \right)
  \Big(
   \varphi(u(t)) + |u(t)|_V^p
  \Big)
\label{At:xdu}
\end{eqnarray}
for a.e. $t \in (0, \infty)$.
Multiplying (\ref{At:xu}) by $(16C_7)^{-1}$ and adding this to
 (\ref{At:xdu}), we derive
\begin{eqnarray*}
  \frac{d}{dt} 
  \Big(
   \varphi(u(t)) + \sigma |u(t)|_V^p
  \Big)
  + \left(
     \frac{\alpha}{32 C_7} 
     - \lambda^{p'} M_2 - \frac{p}{p'} \sigma
    \right)
  \Big(
   \varphi(u(t)) + |u(t)|_V^p
  \Big)
\nonumber \\
 \leq C \left(
	|f|_{V^*}^{p'} + 1
       \right)
\end{eqnarray*}
for a.e. $t \in (0, \infty)$.
Set $\phi(t) := \varphi(u(t)) + \sigma |u(t)|_V^p$. 
If $\lambda > 0$ is so small that
\begin{eqnarray*}
 && 
\beta := 
\frac{\alpha}{32 C_7} 
- \lambda^{p'} M_2 - \frac{p}{p'} \sigma
\geq
\frac{\alpha}{64 C_7} - \lambda^{p'} M_2
> 0,
\end{eqnarray*}
then 
\begin{eqnarray}\label{At:phi-odi}
 && \frac{d \phi}{dt}(t) + \beta \phi(t) \leq F := C \left(
	|f|_{V^*}^{p'} + 1
       \right)
 \quad \mbox{ for a.e. } \ t \in (0, \infty).
\end{eqnarray}
By using standard techniques for differential inequalities, we have
\begin{eqnarray*}
 && \phi(t) \leq \frac{F}{\beta} + \phi(0) e^{-\beta t}
  \quad \mbox{ for all }\ t \in [0, \infty),
\end{eqnarray*}
in particular,
\begin{eqnarray*}
 && \phi(t) \leq \frac{F}{\beta} + 1
  \quad \mbox{ for all } \ t \geq \log (\phi(0) + 1)/\beta.
\end{eqnarray*}
Thus, by putting $R := F/(\sigma\beta) + 1/\sigma$ and $T_0(\cdot) := \log (\cdot + 1)/\beta$,
 we obtain (\ref{At:phi-bdd}).
\end{proof}

The dissipative estimate obtained above implies the point dissipativity
and the eventual boundedness of $\mathcal{G}$.

\begin{lemma}\label{At:L:pd+eb}
Under the same assumptions as in Theorem {\rm \ref{At:T:At}}, the
 following {\rm (i)} and {\rm (ii)} are satisfied.
 \begin{enumerate}
  \item[\rm (i)] $\mathcal{G}$ is point dissipative.
  \item[\rm (ii)] $\mathcal{G}$ is eventually bounded.
 \end{enumerate}
\end{lemma}

\begin{proof}
Let $R \geq 0$ and $T_0(\cdot)$ be the constant and the increasing function
 given by Lemma \ref{At:L:phi-bdd} respectively. Moreover, we write
\begin{equation*}
 B_r := \{v \in X ; \ \varphi(v) + |v|_V^p \leq r\}
  \quad \mbox{ for } \ r > 0.
\end{equation*}

\noindent
{\it Proof of {\rm (i)}.} Put $B := B_R$. 
Let $u \in \mathcal{G}$ and set $\tau := T_0(\varphi(u(0)) + |u(0)|_V^p)$. 
Then by Lemma \ref{At:L:phi-bdd}, we can deduce that $u(t) \in B$ for all $t \geq \tau$.   

\noindent
{\it Proof of {\rm (ii)}.} 
Let $D$ be a bounded set in $X$. Then we can take $R_1 \in (0, \infty)$ such that 
$D \subset B_{R_1}$. Moreover, put $\tau := T_0 (R_1)$ and $B := B_R$. 
Then by Lemma \ref{At:L:phi-bdd}, for any $u \in \mathcal{G}$ with $u(0) \in D \subset B_{R_1}$, 
it follows that $u(t) \in B$ for all $t \geq \tau$.
\end{proof}

Concerning the compactness of $\mathcal{G}$, we have:

\begin{lemma}\label{At:L:cpt}
Under the same assumptions as in Theorem {\rm \ref{At:T:At}},
 $\mathcal{G}$ is compact.
\end{lemma}

\begin{proof}
Since $u_n(0)$ is bounded in $X$, i.e., $|u_n(0)|_V + \varphi(u_n(0)) \leq C$ 
with some constant $C$ independent of $n$, the estimates
 (\ref{GF:e:u+phi})--(\ref{GF:e:du}), (\ref{GF:e:dpsi})--(\ref{GF:e:dphi})
 and the convergences (\ref{GF:c:ut}), (\ref{GF:c:du})--(\ref{GF:c:B-s})
 are established with an arbitrary $T > 0$ as in the proof of (H4) (see
 the proof of Theorem \ref{GF:T:GF}). 
Moreover, we can also verify that there exist a subsequence $(n_k)$ of
 $(n)$ and $u \in C([0,\infty);X)$ such that
\begin{eqnarray}\label{At:phi-conv}
 && \varphi(u_{n_k}(t)) \to \varphi(u(t))
  \quad \mbox{ for all }\ t \in (0, \infty).
\end{eqnarray}
Thus $u_{n_k}(t) \to u(t)$ in $X$ for each positive $t$.
This completes our proof.
\end{proof}

\begin{proof}[Proof of Theorem \ref{At:T:At}]
Thanks to Theorem \ref{P:T:At} and Proposition \ref{P:P:asym-cpt}, 
we can immediately derive the conclusion of Theorem \ref{At:T:At} from Lemmas
 \ref{At:L:phi-bdd}--\ref{At:L:cpt}.
\end{proof}

\begin{remark}\label{At:R:crit}
{\rm
 Theorem \ref{At:T:At} treats a critical case that the perturbation
 $B(u)$ of (CP) almost balances with the leading term $\partial_V \varphi(u)$
 in a sense and the domination will be determined by the parameter
 $\lambda > 0$.
 Indeed, the following concrete PDE problem can be reduced into (CP)
 concerned in Theorem \ref{At:T:At}:
 $$
 u_t - \Delta u = \lambda u \ \mbox{ in } \ \Omega \times (0, \infty),
 \ u = 0 \ \mbox{ on } \ \partial \Omega \times (0, \infty),
 \ u(\cdot, 0) = u_0 \ \mbox{ in } \ \Omega,
 $$
 where $\Omega$ is a bounded domain of $\mathbb R^N$,
 by putting $X_0 = X = H^1_0(\Omega)$, $V = H = V^* = L^2(\Omega)$,
 $\partial_V \psi(u) = u$, $\partial_V \varphi(u) = -\Delta u$ and 
 $B(u) = -u$. Then the unique solution $u = u(x,t)$ of (CP) for some
 data is unbounded as $t \to \infty$ (hence there is no
 bounded absorbing set) for the case $\lambda > \lambda_1$, where
 $\lambda_1 = \lambda_1 (\Omega)$ denotes the principal eigen-value of
 $-\Delta$ with the homogeneous Dirichlet boundary condition, and every
 solution converges to the origin for the case $\lambda < \lambda_1$. 
 Moreover, we also note that the inequality of (S1) holds only for
 $\lambda < \lambda_1$.
}
\end{remark}

Let us next consider the case that the perturbation $B(u)$ is completely
dominated by the leading term $\partial_V \varphi(u)$. Then we can
remove the constraint on the smallness of $\lambda > 0$ in Theorem
\ref{At:T:At}. To do so, we use the following condition 
instead of ${\rm (B1)}_\varepsilon$.
 \begin{enumerate}
  \item[\rm (B1)] $D(\partial_V \varphi) \subset D(B)$.
	       For any $\varepsilon > 0$, there exist constants
	       $r_\varepsilon \in [0, 1)$ and $c_\varepsilon \geq 0$
	       such that
	       \begin{eqnarray*}
		&& |g|_{V^*}^{p'} \leq \varepsilon |\xi|_{V^*}^\sigma 
		 + c_\varepsilon \left\{
				  \left(
				   |\varphi(u)| + |u|_V^p
				  \right)^{r_\varepsilon} + 1 
				 \right\}
		 \quad \mbox{ with } \sigma := \min\{2,p'\}
	       \end{eqnarray*}
	       for all $u \in D(\partial_V \varphi)$, $g \in B(u)$ 
	       and $\xi \in \partial_V \varphi(u)$. 
 \end{enumerate}

More precisely, we have:

\begin{corollary}\label{At:C:At}
Let $p \in (1, \infty)$ be fixed and assume that {\rm (A1), (A2),
 ($\Phi$1)} and {\rm (B1)--(B3)} for any $T > 0$. Moreover, suppose that
 {\rm (S1)} is satisfied with $\bar \lambda = 1$. 
Then the same conclusion as in Theorem {\rm
 \ref{At:T:At}} holds true for any $\lambda \in [0,1]$.
\end{corollary}

\begin{proof}
The smallness of $\lambda$ is used only in the proof of Lemma
 \ref{At:L:phi-bdd}. So it suffices to prove this lemma without the
 smallness by assuming (B1) instead of ${\rm (B1)}_\varepsilon$.
Let $u \in \mathcal{G}$. Then as in Proposition \ref{GE:P:B-bdd}, it
 follows from (B1) that there exists a constant
 $\gamma_\varepsilon > 0$ such that $\gamma_\varepsilon \to 0$ as
 $\varepsilon \to 0$ and
\begin{eqnarray}
|g(t)|_{V^*}^{p'} 
  \leq M_\varepsilon \left\{
	  |f|_{V^*}^{p'} + C_4
	  + \Big(
	     \varphi(u(t)) + |u(t)|_V^p
	 \Big)^{r_\varepsilon}
	  + 1 \right\}
	  + \gamma_\varepsilon \psi \left(u'(t)\right)
\label{GF:B1'}
\end{eqnarray}
for a.e. $t \in [0, \infty)$. Hence as in the proof of Lemma
 \ref{At:L:phi-bdd}, by choosing $\varepsilon > 0$ and $\sigma > 0$ 
small enough,
\begin{eqnarray}
\lefteqn{
 \frac 1{16} \psi(u'(t)) 
+ \frac{d}{dt} \Big(
		\varphi(u(t)) + \sigma |u(t)|_V^p
	       \Big)
}\nonumber \\
 &\leq&
  C \left(
     |f|_{V^*}^{p'} + 1
    \right)
  + \lambda^{p'} M_2 \Big(
		 \varphi(u(t)) + |u(t)|_V^p
	    \Big)^{r_\varepsilon}
  + \frac{p}{p'} \sigma |u(t)|_V^p
\label{At:xdu-B1'}
\end{eqnarray}
for a.e. $t \in (0, \infty)$. Since $r_\varepsilon \in (0, 1)$, 
by Young's inequality, for any $\mu > 0$ we can take a constant $d_\mu > 0$ such that
\begin{eqnarray*}
  &&
   \Big(
    |\varphi(u)| + |u|_V^p
   \Big)^{r_\varepsilon}
   \leq 
   \mu \Big(
	    |\varphi(u)| + |u|_V^p
	   \Big) + d_\mu.
\end{eqnarray*}
Here we particularly choose $\mu > 0$ so small that
\begin{equation*}
\beta := 
 \frac{\alpha}{32 C_7} 
 - \lambda^{p'} \mu M_2 - \frac{p}{p'} \sigma
 \geq 
 \frac{\alpha}{64 C_7}  - \lambda^{p'} \mu M_2
> 0.
\end{equation*}
Then we obtain (\ref{At:phi-odi}) from (\ref{At:xdu-B1'}) and the rest
 of proof runs as before. 
\end{proof}

\section{Applications to generalized Allen-Cahn equations}\label{S:App}

Gurtin~\cite{Gurtin} proposed a generalized Allen-Cahn equation,
which describes the evolution of an order parameter $u = u(x,t)$, of the
form
\begin{equation}\label{Gurtin}
 \rho (u, \nabla u, u_t) u_t 
  = \mathrm{div} \left(
		  \partial_{\bf p} \hat{\psi}(u, \nabla u)
		 \right)
  - \partial_r \hat{\psi}(u, \nabla u) + f,
\end{equation}
where $\rho = \rho(r, {\bf p}, s) \geq 0$ is a kinetic modulus,
$\hat{\psi} = \hat{\psi} (r, {\bf p})$ denotes a free energy density
with its derivatives $\partial_r \hat{\psi}$ and $\partial_{\bf p}
\hat{\psi}$ in $r$ and ${\bf p}$, respectively, and
$f$ is an external microforce, by taking account of a balance law for
microforces as well as constitutive relations consistent with the second
law of thermodynamics. As a simple instance of the free energy density
$\hat{\psi}$, we often take
\begin{equation*}
 \hat{\psi}(r, {\bf p}) = \frac{1}{2} |{\bf p}|^2 + W(r)
\end{equation*}
with a double-well potential $W(r)$ (e.g., $W(r) = (r^2 - 1)^2$). Moreover, if
 $\rho \equiv 1$, then (\ref{Gurtin}) coincides with the usual
 Allen-Cahn equation,
\begin{equation*}
 u_t - \Delta u + \partial_r W(u) = f.
\end{equation*}

In \S \ref{Ss:App:deg}, we particularly deal with a generalized
Allen-Cahn equation of quasilinear type such as
\begin{equation}\label{gAC-1}
\left.
\begin{array}{ll}
 \alpha (u_t) - \Delta_m u + \partial_r W(x, u) \ni
  f & \mbox{ in } \ \Omega \times (0, \infty),\\[1mm]
 u = 0 &  \mbox{ on } \ \partial \Omega \times (0, \infty),\\[1mm]
 u(\cdot,0) = u_0 & \mbox{ in } \  \Omega, 
\end{array}
\right\}
\end{equation}
where $u_0, f: \Omega \to \mathbb{R}$ are given, $\alpha(r) =
|r|^{p-2}r$ with $p \geq 2$, and $\Delta_m$ stands
for the so-called $m$-Laplace operator given by
\begin{equation*}
  \Delta_m u(x) = \nabla \cdot \left(
			       |\nabla u(x)|^{m-2}\nabla u(x)
			      \right),
 \quad 1 < m < \infty.
\end{equation*}
Here we set
\begin{equation*}
\partial_r W(x,r) = \partial j(r) + \lambda g(x,r)
\end{equation*}
with the subdifferential $\partial j$ of a proper, lower semicontinuous convex
function $j : \mathbb{R} \to [0, +\infty]$, $\lambda \in [0,1]$ and a 
Carath\'eodory function $g : \Omega \times \mathbb{R} \to \mathbb{R}$ possibly
non-monotone in the second variable.
Then $\partial_r W$ can be regarded as the derivative (in a proper sense) 
of a potential $W = W(x,r) : \Omega \times
\mathbb{R} \to (-\infty, +\infty]$ given by
\begin{equation}\label{App:W-typ}
W(x,r) := j(r) + \lambda \int^r_0 g(x,\rho) d\rho
\quad \mbox{ for }\ x \in \Omega, \ r \in \mathbb{R}.
\end{equation}
As a typical example of $j$, we take a power function, i.e., $j(r) =
|r|^\sigma$, or the indicator function of the unit ball, i.e., $j(r) =
0$ if $|r| \leq 1$; $j(r) = \infty$ if $|r| > 1$.
Here one can easily find that (\ref{gAC-1}) is reduced into
Gurtin's generalized Allen-Cahn equation (\ref{Gurtin}) with a spatially
inhomogeneous free energy density $\hat{\psi} = \hat{\psi}(x,r,{\bf p})$
by setting
\begin{equation*}
 \hat{\psi}(x, r, {\bf p}) = \frac{1}{m} |{\bf p}|^m + W(x, r)
\quad \mbox{ and } \quad
\rho(r,{\bf p},s) = |s|^{p-2}.
\end{equation*}

In \S \ref{Ss:App:gn} we also consider the following equation as a
further generalized form of semilinear Allen-Cahn equations.
\begin{equation}\label{gAC-2}
\left.
\begin{array}{ll}
 \alpha (u_t) - \Delta u + N(x, u, \nabla u) = f &  
  \mbox{ in } \ \Omega \times (0, \infty),\\[1mm]
 u = 0 &  \mbox{ on } \ \partial \Omega \times (0, \infty),\\[1mm]
 u(\cdot, 0) = u_0 & \mbox{ in } \ \Omega, 
\end{array}
\right\}
\end{equation}
where $N = N(x, r, {\bf p})$ is written as follows
\begin{equation*}
 N(x,r,{\bf p}) = \partial j(r) + \lambda h(x, r, {\bf p})
\end{equation*}
with a Carath\'eodory function $h : \Omega \times \mathbb{R} \times
\mathbb{R}^N \to \mathbb{R}$, and $\partial j$ is single-valued.
We emphasize that this problem may not be written as a (generalized)
gradient system such as (\ref{Gurtin}), since the nonlinear term $N$
depends on the gradient of $u$. However, this problem can be regarded as
a perturbation problem of (\ref{Gurtin}).

We shall apply the preceding abstract theory to (\ref{gAC-1}) and
(\ref{gAC-2}) to prove
the existence of global attractors in Ball's sense.
Throughout this section, let $\Omega$ be a bounded domain of
$\mathbb{R}^N$ with smooth (e.g., $C^2$-class) boundary $\partial \Omega$. 

\subsection{Quasilinear Allen-Cahn equations}\label{Ss:App:deg}

The aim of this subsection is to discuss the existence of global (in
time) solutions for (\ref{gAC-1}) and their asymptotic behavior. 
We assume that the mappings $\alpha : \mathbb{R} \to \mathbb{R}$
and $g : \Omega \times \mathbb{R} \to \mathbb{R}$ satisfy 
\begin{enumerate}
 \item[(a1)] there exists $p \in [2, \infty)$ such that $\alpha (r) =
	     |r|^{p-2}r$ \ for $r \in \mathbb{R}$;
 \item[(a2)] $j = j(r)$ is a proper lower semicontinuous convex
	     function from $\mathbb{R}$ into $[0, \infty]$ such
	     that $\partial j(0) \ni 0$; there
	     exist constants $\sigma > 1$, $C_8 > 0$ and $C_9 \geq 0$ such that
	     \begin{equation*}
	      C_8 |r|^\sigma \leq j(r) + C_9
	       \quad \mbox{ for all } \ r \in \mathbb{R};
	     \end{equation*}
 \item[(a3)] $g = g(x,r)$ is a Carath\'eodory function, i.e., measurable
	     in $x$ and continuous in $r$. Moreover, there exist $q \geq
	     1 + 1/p'$, $C_{10} \geq 0$ and $a_1 \in L^1(\Omega)$ such that
	     \begin{equation*}
	      |g(x,r)|^{p'} \leq C_{10} |r|^{p'(q-1)} + a_1(x)
	     \end{equation*}
	     for a.e. $x \in \Omega$ and all $r \in \mathbb{R}$.
\end{enumerate}
Our basic assumption on the exponent $p$ is as follows:
\begin{equation}\label{App:b-hyp}
 2 \leq p < \max \{m^*, \sigma\},
\end{equation}
where $m^*$ denotes the Sobolev critical exponent, that is,
\begin{equation*}
 m^* := \frac{Nm}{(N - m)_+}
  = \left\{
	 \begin{array}{ll}
	  \dfrac{Nm}{N - m} \ &\mbox{ if } \ m < N,\\
	  \infty \ &\mbox{ if } \ m \geq N.
	 \end{array}
	\right.
\end{equation*}

In order to reduce (\ref{gAC-1}) into an abstract Cauchy problem, 
we set $V = L^p(\Omega)$, $V^* = L^{p'}(\Omega)$ and $H = L^2(\Omega)$.
Then (\ref{ps}) follows.
Moreover, define the functional $\psi \colon V \to [0, \infty)$ by
\begin{equation}
 \psi(u) := \frac{1}{p} \int_\Omega |u(x)|^p dx \quad \mbox{ for } \ u \in V.
  \label{App:psi}
\end{equation}
Then $D(\psi) = D(\partial_V \psi) = V$ and $\partial_V \psi(u)$ coincides with $\alpha(u(\cdot))$ in $V^*$.
Define $\varphi_1, \varphi_2 \colon V \to [0, \infty]$ by
\begin{equation*}
 \varphi_1(u) := \left\{
		  \begin{array}{ll}
		   \dfrac{1}{m} \displaystyle \int_\Omega |\nabla u(x)|^m
		    dx 
		    & \mbox{ if } \ u \in W^{1,m}_0(\Omega), \\[2mm]
		    \infty \quad & \mbox{ otherwise}
		  \end{array}
		 \right.
\end{equation*}
and
\begin{equation*}
 \varphi_2(u) := \left\{
		  \begin{array}{ll}
		   \displaystyle \int_\Omega j(u(x)) dx 
		    & \mbox{ if } \ j(u(\cdot)) \in L^1 (\Omega), \\[2mm]
		    \infty \quad & \mbox{ otherwise. }
		  \end{array}
		 \right.
\end{equation*}
It is easily seen that $\partial_V \varphi_1 (u) = - \Delta_m u$ with
$u|_{\partial \Omega} = 0$ and $\partial_V \varphi_2 (u) = \partial
j(u(\cdot))$.
We further set $\varphi: V \to [0, \infty]$ by
\begin{equation}\label{App:phi}
 \varphi(u) := \varphi_1(u) + \varphi_2(u)
\end{equation}
with the effective domain $D(\varphi) = D(\varphi_1) \cap D(\varphi_2) =
\{u \in W^{1,m}_0(\Omega); \ j(u(\cdot)) \in L^1(\Omega)\}$.
Then since $\partial j(0) \ni 0$, as in Corollary 16 of~\cite{HB3}, we
observe
\begin{equation*}
 \varphi_1(J^H_\mu u) \leq \varphi_1(u)
  \quad \mbox{ for } \ u \in D(\varphi_1) \ \mbox{ and } \ \mu > 0,
\end{equation*}
where $J^H_\mu$ denotes the resolvent of the subdifferential for the
extension of $\varphi_2$ onto $H$ by infinity. Then we can verify that
$\partial_V \varphi_1 + \partial_V \varphi_2$ is maximal monotone;
hence, $\partial_V \varphi(u) = \partial_V \varphi_1(u) + \partial_V
\varphi_2(u)$ coincides with $- \Delta_m u + \partial j(u(\cdot))$ equipped
with $u|_{\partial \Omega} = 0$ in $V^*$ (see~\cite{G15}). Furthermore, we define the mapping
\begin{equation*}
B : V \to V^*; \ u \mapsto g(\cdot, u(\cdot)).
\end{equation*}
Then (\ref{gAC-1}) is rewritten into (CP) with the functionals $\psi, \varphi$
and the mapping $B$ defined above. 

In the rest of this paper, a function $u: \Omega \times
(0,\infty) \to \mathbb{R}$ is said to be an {\it $L^p$-solution} of
(\ref{gAC-1}) on $[0,\infty)$, if $u$ is a strong solution on $[0, \infty)$ of
(CP) with $\psi, \varphi, B$ defined above.
Applying the preceding abstract theory to (\ref{gAC-1}), we can prove the
following theorems.
\begin{theorem}\label{App:T:GEGF}
Let $\Omega$ be a bounded domain of $\mathbb{R}^N$ with smooth boundary
 $\partial \Omega$. Suppose that {\rm (a1)--(a3)} and {\rm
 (\ref{App:b-hyp})} are satisfied, and 
\begin{equation}\label{App:hypo}
p'(q - 1) \leq \max \left\{
		 m, p
		 \right\}
\quad \mbox{or} \quad
 \left\{
  \begin{array}{ll}
   p'(q - 1) \leq \sigma & \mbox{ if } \ \sigma < m^*,\\
   p'(q - 1) < \sigma & \mbox{ if } \ \sigma \geq m^*.
  \end{array}
 \right.
\end{equation}
Then for any $\lambda \in [0,1]$, $f \in L^{p'}(\Omega)$, $u_0 \in
 X := \{ v \in W^{1,m}_0(\Omega);\ j(v(\cdot)) \in
 L^1(\Omega)\}$, the
 initial-boundary value problem {\rm (\ref{gAC-1})} admits at least one
 $L^p$-solution on $[0, \infty)$. Moreover, the set
\begin{eqnarray*}
 \mathcal{G} := \{
		 u : [0, \infty) \to L^p(\Omega) ;\
		 u \mbox{ is an $L^p$-solution of {\rm (\ref{gAC-1})} on } [0,
		 \infty)\quad
		 \\
		 \mbox{ with some } u_0 \in X
		\}
\end{eqnarray*}
forms a generalized semiflow in $X$ with the distance $d_X(\cdot, \cdot)$ given by
\begin{equation*}
 d_X(u, v) := \left|
	       \varphi(u) - \varphi(v)
	      \right|
 + |u - v|_{L^p(\Omega)}
 \quad \mbox{ for } \ u \in X,
\end{equation*}
where $\varphi$ is given by {\rm (\ref{App:phi})}.
\end{theorem}

\begin{remark}\label{App:R:hypo}
 {\rm
We can derive from the assumptions (\ref{App:b-hyp}) and (\ref{App:hypo}) that
\begin{enumerate}
 \item[(i)] $p'(q - 1) \leq \max\{m, p, \sigma\}$,
 \item[(ii)] $X_0 := W^{1,m}_0(\Omega) \cap L^\sigma(\Omega)$ is compactly
	     embedded both in $L^p(\Omega)$ and in $L^{p'(q - 1)}(\Omega)$.
\end{enumerate}
Indeed, (i) follows immediately from (\ref{App:hypo}). Moreover,
 applying the Rellich-Kondrachov compactness theorem and interpolation inequalities
 for Lebesgue spaces, we assure by (\ref{App:b-hyp}) that $X_0$ is
 compactly embedded in $L^p(\Omega)$. Furthermore, (ii) is similarly
 derived from (\ref{App:b-hyp}) and (\ref{App:hypo}). 
These two facts play a crucial role in our proofs and we can replace
 (\ref{App:b-hyp}) and (\ref{App:hypo}) by (i) and (ii) in Theorems
 \ref{App:T:GEGF} and \ref{App:T:At}.
}
\end{remark}

\begin{proof}
Conditions (A1) and (A2) follow immediately from (a1) and
 (\ref{App:psi}). Set $X_0 := W^{1,m}_0(\Omega) \cap L^\sigma(\Omega)$
 with the norm 
$|\cdot|_{X_0} := (|\nabla \cdot|_{L^m(\Omega)}^2 + |\cdot|_{L^\sigma(\Omega)}^2)^{1/2}$. 
Then we get, by (a2)
\begin{equation*}
 |u|_{X_0}^{\min\{m, \sigma\}} \leq C \left(
				       \varphi(u) + 1
				      \right)
 \quad \mbox{ for all } \  u \in D(\varphi)
\end{equation*}
with some constant $C > 0$. Moreover, by (ii) of Remark
 \ref{App:R:hypo}, $X_0$ is compactly embedded in $V$.
Hence ($\Phi$1) holds. As for ${\rm (B1)}_\varepsilon$, we
 note by (a3) that
\begin{equation}\label{App:B-bdd}
 |B(u)|_{V^*}^{p'} 
  \leq 
   C_{10} \int_\Omega |u(x)|^{p'(q-1)} dx + \int_\Omega a_1(x) dx
   \quad \mbox{ for } \ u \in L^{p'(q-1)}(\Omega).
\end{equation}
Here by (ii) of Remark \ref{App:R:hypo}, we have 
$D(\varphi) \subset X_0 \subset L^{p'(q-1)}(\Omega) \subset D(B)$, and
 moreover, by (a2)
\begin{equation*}
\int_\Omega |u(x)|^{p'(q-1)} dx
  \leq C \left(
	  \varphi(u) + |u|_V^p + 1 
	 \right)
  \quad \mbox{ for } \ u \in D(\varphi),
\end{equation*}
which together with (\ref{App:B-bdd}) implies ${\rm (B1)}_\varepsilon$. 
Furthermore, by (\ref{App:B-bdd}), the mapping $u \mapsto B(u(\cdot))$
 becomes continuous from $L^{p'(q-1)}(0,T;L^{p'(q-1)}(\Omega))$ into $L^{p'}(0,T;V^*)$
 for any $T > 0$ (see Theorem 1.43 of~\cite{Roubicek}), in particular, (B3)
 holds. Recall that $X_0$ is compactly embedded in
 $L^{p'(q-1)}(\Omega)$ by (ii) of Remark \ref{App:R:hypo}. Hence we also obtain (B2),
 since the sequence $(u_n)$ of (B2) becomes precompact in
 $C([0,S];L^{p'(q-1)}(\Omega))$ by Aubin-Lions's compactness theorem.
Finally, applying Theorems \ref{GE:T:GE} and \ref{GF:T:GF}, we
 can obtain our desired conclusion.
\end{proof}

\begin{theorem}\label{App:T:At}
Let $\Omega$ be a bounded domain of $\mathbb{R}^N$ with 
smooth boundary $\partial \Omega$ and let $f \in L^{p'}(\Omega)$. 
In addition to {\rm (a1)--(a3)}, {\rm (\ref{App:b-hyp})} and {\rm
 (\ref{App:hypo})}, suppose that
\begin{equation}
p \leq \max\{m,\sigma\}
\quad \mbox{and} \quad 
\left\{
\begin{array}{ll}
\lambda \mbox{ is arbitrary } &\mbox{if } \ p'(q - 1) < \max\{m, \sigma\},\\
\lambda \mbox{ is small } &\mbox{otherwise.}
\end{array}
\right.
\label{App:hypo-2}
\end{equation}
Then the generalized semiflow $\mathcal{G}$ has a global
 attractor $\mathcal{A}$ in $X$.
\end{theorem}

\begin{proof}
We can derive from (a3), (i) of Remark \ref{App:R:hypo} and (\ref{App:hypo-2}) that
\begin{eqnarray*}
\lefteqn{
 \left\langle \partial_V \varphi(u) + \lambda B(u) , u \right\rangle
}\\
  &\geq& \int_\Omega |\nabla u(x)|^m dx 
   + \int_\Omega \Big( j(u(x)) -  j(0) \Big) dx
   + \lambda \int_\Omega g(x,u(x)) u(x) dx
   \\
 &\geq& 
  \int_\Omega |\nabla u(x)|^m dx 
  + C_8
  \int_\Omega |u(x)|^\sigma dx - \left( C_9 + |j(0)| \right) |\Omega|
  \\
 && - \lambda \int_\Omega \left(
			 C_{10} |u(x)|^{p'(q-1)} + a_1(x)
			\right)^{1/p'} |u(x)| dx
  \\
 &\geq&
  \frac{1}{2} \int_\Omega |\nabla u(x)|^m dx 
  + \frac{C_8}{2} \int_\Omega |u(x)|^\sigma dx - M_\lambda
\end{eqnarray*}
with some constant $M_\lambda$ depending on $\lambda$,
by assuming that $\lambda > 0$ is small enough only in case 
$\max\{m, \sigma\} = q$. 
Hence (S1) follows from the fact that $p \leq \max\{m, \sigma\}$.
Consequently, Theorem \ref{At:T:At} proves the existence of a global
 attractor with a sufficiently small $\lambda > 0$. 

In addition, assuming $p'(q - 1) < \max\{m, \sigma\}$ 
(then, it is not true that $\max\{m, \sigma\} = q$), we observe
\begin{equation*}
 \int_\Omega |u(x)|^{p'(q-1)} dx
  \leq C \left(
	  \varphi(u) + |u|_V^p + 1 
	 \right)^r
  \quad \mbox{ for } \ u \in D(\varphi)
\end{equation*}
with some $r \in (0, 1)$. Hence (B1) follows.
Therefore Corollary \ref{At:C:At} ensures the existence of a global
 attractor for any $\lambda \in [0,1]$.
\end{proof}

\begin{remark}
 {\rm
\begin{enumerate}
 \item[(i)] Our arguments described above are still valid even if (a1)
	    is relaxed into 
	    \begin{enumerate}
	     \item[$({\rm a1})'$] there exist $p \in [2, \infty)$ and a lower semicontinuous convex
			  function $A : \mathbb{R} \to [0,\infty)$ such
			  that $\partial A = \alpha$ and
			  \begin{equation*}
			   C_{11} (|r|^p - 1) \leq A(r)
			    \ \mbox{and} \
			    |\eta|^{p'} \leq C_{12}(A(r) + 1) 
			   \ \mbox{ for all } [r, \eta] \in \alpha
			  \end{equation*}
			  with some constants $C_{11}, C_{12} \geq 0$.
	    \end{enumerate}
	    Furthermore, we can also treat the case where $-\Delta_m u(x,t)$
	    is replaced by $- \mathrm{div}~ {\bf a}(x, \nabla u(x,t))$ and
	    ${\bf a} = {\bf a}(x,{\bf p})$ is a function from
	    $\Omega \times \mathbb{R}^N$ into $\mathbb{R^N}$ such that
	    ${\bf a}$ is continuous and maximal monotone in ${\bf p}$
	    and measurable in $x$, by imposing appropriate growth
	    conditions on ${\bf a}$ in ${\bf p}$ (see~\cite{G11} for
	    more details).
\item[(ii)] If we have a boundedness condition of $\partial j$, we could
	    weaken the assumptions on $q$, the growth order of the
	    perturbation. Let us consider a simple case that $j(r) = (1/\sigma) |r|^\sigma$ 
	    with $\sigma > 1 + 1/p'$.
	    Then $\partial j(r) = |r|^{\sigma - 2}r$, which implies
	    $|\partial j(r)| = |r|^{\sigma - 1}$. Here we note that
	    \begin{equation*}
	     |u|_{L^{p'(\sigma-1)}(\Omega)}^{\sigma - 1}
	      =
	      |\partial j(u(\cdot))|_{L^{p'}(\Omega)}
	      \leq |-\Delta_m u + \partial j(u(\cdot))|_{L^{p'}(\Omega)}
	    \end{equation*}
	    for $u \in W^{1,m}_0(\Omega) \cap L^{p'(q - 1)}(\Omega)$
	    (see~\cite{G15}). Hence in case $q < \sigma$, for any
	    $\varepsilon > 0$, we can choose $C_\varepsilon > 0$ such that
	    \begin{equation*}
	     |u|_{L^{p'(q-1)}(\Omega)}^{p'(q-1)}
	      \leq \varepsilon
	      |u|_{L^{p'(\sigma-1)}(\Omega)}^{p'(\sigma-1)} +
	      C_\varepsilon
	      \leq \varepsilon |\partial_V \varphi(u)|_{V^*}^{p'} + C_\varepsilon.
	    \end{equation*}
	    Combining this fact with (\ref{App:B-bdd}), we obtain (B1).
	    Moreover, under the same assumption,
	    $D := W^{1,m}_0(\Omega) \cap L^\sigma(\Omega) \cap L^{p'(\sigma - 1)}(\Omega)$ is
	    compactly embedded in $L^{p'(q-1)}(\Omega)$; hence the
	    sequence $(g_n)$ of (B2) becomes precompact in $L^{p'}(0,S;V^*)$,
	    since $D(\partial_V \varphi) \subset D$. Thus we can
	    prove Theorem \ref{App:T:GEGF} with $q < \sigma$ instead
	    of (\ref{App:hypo}). Furthermore, as in Theorem
	    \ref{App:T:At}, we can also verify the existence of a global
	    attractor for any $\lambda \in [0,1]$ under $q < \sigma$
	    instead of (\ref{App:hypo}) by assuming a structure
	    condition such as
	    \begin{equation*}
	     g(x,r) r \geq - C \left( |r|^q + 1 \right)
	      \quad \mbox{ for a.e. } \ x \in \Omega 
	      \ \mbox{ and } \ r \in \mathbb{R}. 
	    \end{equation*}
\end{enumerate} 
}
\end{remark}

\subsection{Semilinear Allen-Cahn equations with perturbations
  involving gradients}\label{Ss:App:gn}

We next treat (\ref{gAC-2}). Let us assume
\begin{enumerate}
 \item[${\rm (a2)}'$] $\partial j$ is single-valued, and (a2) holds;
 \item[(a4)\;] The function $h = h(x,r, {\bf p})$ is a Carath\'eodory
	     function in $\Omega \times \mathbb{R} \times \mathbb{R}^N$ 
	     (i.e., measurable in $x$ and continuous in $(r, {\bf p})$). 
	     Moreover, there exist constants $q_1, q_2 \geq 1 + 1/p'$, 
	     $C_{13} \geq 0$ and  
	     $a_2 \in L^1(\Omega)$ such that 
	     \begin{eqnarray*}
	      && |h(x,r,{\bf p})|^{p'} \leq C_{13} (|r|^{p'(q_1 - 1)} 
	       + |{\bf p}|^{p'(q_2 - 1)}) + a_2(x)
	     \end{eqnarray*}
	     for a.e. $x \in \Omega$ and all $(r,{\bf p}) \in \mathbb{R} \times \mathbb{R}^N$.
\end{enumerate}
Furthermore, we use the following basic assumption on the exponent $p$,
\begin{equation}
2 \leq p < \max\{2^*, \sigma\}.
\label{App:b-hyp2}
\end{equation}

Equation (\ref{gAC-2}) has no longer any gradient structure, because the
nonlinear term $N$ depends on $\nabla u$ also. However, (\ref{gAC-2})
can be transcribed into a (non-monotone) perturbation problem for a
doubly nonlinear gradient system in the form (\ref{ee}). Indeed, we set
$V = L^p(\Omega)$, $H = L^2(\Omega)$ and $V^* := L^{p'}(\Omega)$ as well
as functionals $\psi$ and $\varphi$ as in (\ref{App:psi}) and
(\ref{App:phi}) with $m = 2$ respectively. Moreover, we put
\begin{equation*}
 B(u) := h(\cdot, u(\cdot), \nabla u(\cdot)) \quad \mbox{ for }\
  u \in V.
\end{equation*}
Then (\ref{gAC-2}) is reduced into (CP). Hence $L^p$-solutions of
(\ref{gAC-2}) are also defined as in the case of (\ref{gAC-1}). 
Now, exploiting Theorems \ref{GE:T:GE} and \ref{GF:T:GF}, we can verify:
\begin{theorem}\label{App:T:GEGF2}
Let $\Omega$ be a bounded domain of $\mathbb{R}^N$ with smooth boundary
 $\partial \Omega$. Suppose that {\rm (a1), $(a2)'$, (a4)} and {\rm
 (\ref{App:b-hyp2})} are satisfied. In addition, assume that
\begin{equation}\label{App:hyp3-1}
\mbox{either } q_1 \leq p
\mbox{ or }
\left\{
 \begin{array}{ll}
  p'(q_1 - 1) \leq \sigma & 
   \mbox{ if }\ \sigma < \max\{ 2^*, (p')^\dagger \},\\
  p' (q_1 - 1) < \sigma & 
   \mbox{ if }\ \sigma \geq \max\{ 2^*, (p')^\dagger \}\
 \end{array}
\right.
\mbox{holds,}
\end{equation}
and
\begin{equation}\label{App:hyp3-2}
\left\{
 \begin{array}{ll}
  p'(q_2 - 1) \leq 2 & 
   \mbox{ if }\ 2 < (p')^*,\\
  p'(q_2 - 1) < 2 & 
   \mbox{ if }\ 2 \geq (p')^*,
 \end{array}
\right.
\end{equation}
where $(p')^\dagger := Np'/(N - 2p')_+$ and $(p')^* := Np'/(N - p')_+$.
Then for any $u_0 \in X := \{v \in H^1_0(\Omega) ; \ j(v(\cdot)) \in
 L^1(\Omega) \}$, $f \in L^{p'}(\Omega)$ 
and $\lambda \in [0,1]$, the initial-boundary value
 problem {\rm (\ref{gAC-2})} admits at least one $L^p$-solution on $[0, \infty)$. 
Moreover, the set
\begin{eqnarray*}
 \mathcal{G} := \{
		 u : [0, \infty) \to L^p(\Omega) ;\
		 u \mbox{ is an $L^p$-solution of {\rm (\ref{gAC-2})} on } [0,
		 \infty)\quad
		 \\
		 \mbox{ with some } u_0 \in X
		\}
\end{eqnarray*}
forms a generalized semiflow in $X$ with the distance $d_X(\cdot, \cdot)$ given by
\begin{equation*}
d_X(u, v) 
 := 
\left|
  \varphi(u) - \varphi(v)
\right| 
+ |u - v|_{L^p(\Omega)}
 \quad \mbox{ for } \ u \in X,
\end{equation*}
where $\varphi$ is given by {\rm (\ref{App:phi})} with $m = 2$.
\end{theorem}

\begin{remark}\label{App:R:h1}
 {\rm
As in Remark \ref{App:R:hypo}, assumptions (\ref{App:b-hyp2}),
 (\ref{App:hyp3-1}) and (\ref{App:hyp3-2}) yield in particular that
\begin{enumerate}
 \item[(i)] $p'(q_1 - 1) \leq \max\{p, \sigma\}$ and $p'(q_2 - 1) \leq 2$,
 \item[(ii)] $H^1_0(\Omega) \cap L^\sigma(\Omega)$ is
 compactly embedded in $L^p(\Omega)$, and furthermore, 
	     $W^{2,p'}(\Omega) \cap H^1_0(\Omega) \cap L^\sigma(\Omega)$
	     is compactly embedded in 
	     $L^{p'(q_1 - 1)}(\Omega) \cap W^{1, p'(q_2 - 1)}(\Omega)$.
\end{enumerate}
We can prove Theorems \ref{App:T:GEGF2} and \ref{App:T:At} by using
 these two facts instead of (\ref{App:b-hyp2}), (\ref{App:hyp3-1}) and
 (\ref{App:hyp3-2}). 
}
\end{remark}

\begin{proof}
Due to Theorems \ref{GE:T:GE} and \ref{GF:T:GF}, it suffices to check
 (A1), (A2), ($\Phi$1), ${\rm (B1)}_\varepsilon$--(B3).
As in the proof of Theorem \ref{App:T:GEGF}, we can derive (A1) and
 (A2) from (a1) and (\ref{App:psi}). 
Put $X_0 = H^1_0(\Omega) \cap L^\sigma(\Omega)$. Then
 ($\Phi$1) follows from the fact that $p < \max \{2^*, \sigma\}$ by
 (\ref{App:b-hyp2}). Here we observe by (i) of Remark \ref{App:R:h1} that 
$D(\varphi) \subset H^1_0(\Omega) \cap L^\sigma(\Omega) 
\subset L^{p'(q_1 - 1)}(\Omega) \cap W^{1,p'(q_2 - 1)}(\Omega) \subset D(B)$,
and moreover,
\begin{equation}\label{App:q1q2}
 |u|_{L^{p'(q_1 - 1)}(\Omega)}^{p'(q_1 - 1)} 
  +  |\nabla u|_{L^{p'(q_2 - 1)}(\Omega)}^{p'(q_2 - 1)}
  \leq C \left(
	  \varphi(u) + |u|_V^p + 1
	 \right)
\end{equation}
for all $u \in D(\varphi)$. On the other hand, (a4) yields
\begin{equation}\label{App:B-bdd2}
 |B(u)|_{V^*}^{p'} 
  \leq 
  C_{13} 
  \left(
   |u|_{L^{p'(q_1 - 1)}(\Omega)}^{p'(q_1-1)} + |\nabla u|_{L^{p'(q_2 - 1)}(\Omega)}^{p'(q_2-1)}
  \right)
   + |a_2|_{L^1(\Omega)}
\end{equation}
for $u \in L^{p'(q_1 - 1)}(\Omega) \cap W^{1,p'(q_2-1)}(\Omega)$. 
Thus (\ref{App:q1q2}) implies ${\rm (B1)}_\varepsilon$ with $\varepsilon = 0$.
By Theorem 1.43 of~\cite{Roubicek} and (\ref{App:B-bdd2}), the function
 $\mathcal{B} : u \mapsto B(u(\cdot))$ is continuous from 
\begin{equation*}
 \mathcal{X}_T := L^{p'(q_1 - 1)} (0, T; L^{p'(q_1 - 1)} (\Omega))
 \cap L^{p'(q_2 - 1)}(0, T; W^{1,p'(q_2 - 1)} (\Omega)) 
\end{equation*}
into $L^{p'}(0,T;V^*)$, which particularly implies (B3).

Finally, let us show (B2).
Due to Theorem 9.15 and Lemma 9.17 of~\cite{GT}, we can derive 
$D(\partial_V \varphi_1) = W^{2,p'}(\Omega) \cap H^1_0(\Omega)$ and
\begin{equation*}
|u|_{W^{2,p'}(\Omega)} \leq C |\partial_V \varphi_1(u)|_{V^*} 
  \quad \mbox{ for all } \ u \in D(\partial_V \varphi_1).
\end{equation*}
Hence as in~\cite{G15}, we obtain
\begin{eqnarray}\label{App:CZineq}
 && |u|_{W^{2,p'}(\Omega)} \leq C |\partial_V \varphi(u)|_{V^*} 
  \quad \mbox{ for all } \ u \in D(\partial_V \varphi).
\end{eqnarray}
We note that 
$D(\partial_V \varphi) \subset W^{2,p'}(\Omega) \cap H^1_0(\Omega) \cap 
L^{\sigma}(\Omega)$. 
By (ii) of Remark \ref{App:R:h1}, it follows that $D(\partial_V \varphi)$ 
is compactly embedded in $L^{p'(q_1 - 1)}(\Omega) \cap W^{1,p'(q_2 - 1)}(\Omega)$. 
Hence by (\ref{App:CZineq}) and the
 Aubin-Lions-type compactness theorem (see, e.g., Theorem
 5 of~\cite{Simon}), the sequence $(u_n)$ of (B2) is precompact in 
\begin{equation*}
L^{p'} (0, S; L^{p'(q_1 - 1)} (\Omega)) 
 \cap L^{p'}(0, S; W^{1,p'(q_2 - 1)} (\Omega)). 
\end{equation*}
Moreover, from (\ref{App:q1q2}) and the fact that
\begin{equation*}
 \sup_{t \in [0, S]} \Big(
		      \varphi(u_n(t)) + |u_n(t)|_V
		     \Big) \leq C,
\end{equation*}
the sequence $(u_n)$ is bounded
 in $L^\infty(0,S; L^{p'(q_1 - 1)}(\Omega) \cap W^{1,p'(q_2 - 1)}(\Omega))$. Hence 
we can take a subsequence $(n')$ of $(n)$ such that
\begin{equation*}
u_{n'} \to u \mbox{ strongly in } \mathcal{X}_S.
\end{equation*}
Therefore we infer that 
$\mathcal{B} (u_{n'}) \to \mathcal{B} (u)$ strongly in $L^{p'}(0,T;V^*)$ from
 the continuity of $\mathcal{B}$. Thus (B2) holds.
\end{proof}

Furthermore, the existence of global attractors for $\mathcal{G}$ is
also proved.
\begin{theorem}\label{App:T:At2}
Let $\Omega$ be a bounded domain of $\mathbb{R}^N$ with smooth boundary 
$\partial \Omega$ and let $f \in L^{p'}(\Omega)$. 
In addition to {\rm (a1), $(a2)'$, (a4), (\ref{App:b-hyp2}),
 (\ref{App:hyp3-1})} and {\rm (\ref{App:hyp3-2})}, we suppose that 
\begin{equation*}
p \leq \sigma 
\end{equation*}
and 
\begin{eqnarray}\label{App:hyp3}
\left\{
\begin{array}{ll}
\lambda \mbox{ is arbitrary }  
 &\mbox{ if } \  p'(q_1 - 1) < \sigma \ \mbox{ and } \ p'(q_2 - 1) < 2,
\\
\lambda \mbox{ is small}  &\mbox{ otherwise.}
\end{array}
\right.
\end{eqnarray}
Then the generalized semiflow $\mathcal{G}$ has a global attractor
 $\mathcal{A}$ in $X$.
\end{theorem}

\begin{proof}
From the assumptions, it follows from $(a2)'$ and (a4) that
\begin{eqnarray*}
\lefteqn{
\langle \partial_V \varphi(u) + \lambda B(u), u \rangle
}\\
&\geq& 
\int_\Omega |\nabla u(x)|^2 dx 
+ \int_\Omega \Big( j(u(x)) - j(0) \Big) dx
\\
&& - \lambda \int_\Omega |h(x,u(x),\nabla u(x))| |u(x)| dx
\\
&\geq&
\int_\Omega |\nabla u(x)|^2 dx + C_8 \int_\Omega |u(x)|^\sigma dx 
- \left( C_9 + |j(0)| \right) |\Omega|
\\
&& - \lambda \int_\Omega 
 \left(
  C_{13} |u(x)|^{p'(q_1 - 1)} + C_{13} |\nabla u(x)|^{p'(q_2 - 1)} + a_2(x)
 \right)^{1/p'} |u(x)| dx
 \\
&\geq& 
  \frac{1}{2} \int_\Omega |\nabla u(x)|^2 dx 
  + \frac{C_8}{2} \int_\Omega |u(x)|^\sigma dx
  - M_\lambda
\end{eqnarray*}
for $u \in D(\partial_V \varphi)$. From (i) of Remark \ref{App:R:h1} and
 (\ref{App:hyp3}), the last inequality follows, by assuming that
 $\lambda > 0$ is sufficiently small only in case
\begin{equation}\label{App:sl:small}
\mbox{either } \ q_1 = \sigma \ \mbox{ or } \ q_2 = \frac{2}{\sigma'} + 1.
\end{equation}
Thus (S1) is satisfied, since $p \leq \sigma$.
 
In addition, assuming
\begin{equation}\label{App:sl:any}
p'(q_1 - 1) < \sigma \quad \mbox{ and } \quad p'(q_2 - 1) < 2
\end{equation}
(then, (\ref{App:sl:small}) is not true), we have
\begin{equation*}
  |u|_{L^{p'(q_1 - 1)}(\Omega)}^{p'(q_1 - 1)} 
  +  |\nabla u|_{L^{p'(q_2 - 1)}(\Omega)}^{p'(q_2 - 1)}
  \leq C \left(
	  \varphi(u) + |u|_V^p + 1
	 \right)^r
\end{equation*}
with some $r \in (0, 1)$. Hence (B1) follows from (\ref{App:B-bdd2}).

Therefore Theorem \ref{At:T:At} (respectively, Corollary \ref{At:C:At}) 
ensures that $\mathcal{G}$ admits a global attractor in $X$ for
 sufficiently small (respectively, arbitrary) $\lambda \in [0,1]$
 (respectively, if (\ref{App:sl:any}) holds).
\end{proof}

%% The Appendices part is started with the command \appendix;
%% appendix sections are then done as normal sections
%% \appendix

%% \section{}
%% \label{}

\noindent
{\bf Acknowledgements.}
The author expresses his sincere gratitude to Professors Tomomi Kojo and
 Tomoyuki Idogawa for fruitful discussions. Moreover, the
author also thanks the anonymous referee for careful reading of this
manuscript and for giving useful comments.

\end{document}